\documentstyle[12pt]{article}

\input amssymb.sty

\newcommand{\Two}[4]{\left(\begin{array}{cc}#1&#2\\#3&#4\end{array}\right)}
\newcommand{\Three}[9]{\left(\begin{array}{ccc}#1&#2&#3\\
#4&#5&#6\\#7&#8&#9\end{array}\right)}

\newcommand{\ColFive}[5]{\left(\begin{array}{c}#1\\#2\\#3\\#4\\#5
\end{array}\right)}
\newcommand{\ColF}[5]{\begin{array}{c}#1\\#2\\#3\\#4\\#5\end{array}}
\newcommand{\CFour}[4]{\begin{array}{c}#1\\#2\\#3\\#4\end{array}}

\newcommand{\dia}{\mbox{diag}}
\newcommand{\ant}{\mbox{antidiag}}
\newcommand{\ovr}{\overline}
\newcommand{\np}{\newline\indent}
\newcommand{\ep}{\epsilon}
\newcommand{\al}{\alpha}
\newcommand{\th}{\theta}
\newcommand{\lam}{\lambda}
\newcommand{\pii}{\mbox{\boldmath $\pi$}}
\newcommand{\cir}{\biggl(1-\frac{1}{q}\biggr)}
\newcommand{\EQ}{\begin{equation}}
\newcommand{\EN}{\end{equation}}
\newcommand{\bksl}{\backslash}
\newcommand{\Ch}{C_c(N'\backslash H/K'\!,\psi_{N'})}
\begin{document}

\title{Relation of orbital integrals on $SO(5)$ and $PGL(2)$.}

\author{Dmitrii Zinoviev
\footnote{Department of Mathematics, The Ohio State University,
231 W. 18th Ave., Columbus, OH 43210.}}

\date{}
\maketitle

\begin{abstract} We relate the ``Fourier'' orbital integrals of corresponding
spherical functions on the $p$-adic groups $SO(5)$ and $PGL(2)$.
The correspondence is defined by a ``lifting'' of representations
of these groups. This is a local ``fundamental lemma'' needed
to compare the geometric sides of the global Fourier summation formulae
(or relative trace formulae) on these two groups.
This comparison leads to conclusions about a well known lifting of
representations from $PGL(2)$ to $PGSp(4)$. This lifting produces counter
examples to the Ramanujan conjecture.
\end{abstract}

\vskip0pt plus1in
\penalty-1000
\vskip0pt minus1in

{\bf Introduction.}
Let $G$ be the special orthogonal group, defined
over a local field $F$, by an anti-diagonal form, with upper
triangular minimal parabolic subgroup.
An explicit definition is given in Section 0, where $G$ is denoted by $SO(3,2)$.
Let $C_\th$ ($\th\in F^\times$) be a subgroup of $G$ isomorphic
to the special orthogonal groups $SO(2,2)$ or $SO(3,1)$ (see Section 0).
Denote by $P$ the maximal upper triangular parabolic
subgroup of $G$ with abelian unipotent radical $N$.
For any spherical function $f\in C_c(K\bksl G/K)$, $K$ being
the standard maximal compact subgroup of $G$, consider the
Fourier orbital integral
$$
\int_N\int_{C_\th} f(ngh)\psi_N(n)dhdn,
$$
where $\psi_N$ is a certain character on $N$ depending on a fixed
character $\psi$ of $F$ with conductor $R$ (integers of $F$).
The $N$-$C_\th$-orbits of maximal dimension are of the form
$Na_\al\gamma_0 C_\th$, where $a_\al$ is the diagonal matrix
$\dia(\al,1,1,1,\al^{-1})$ and
$\gamma_0$ is defined below. We are interested in $g$ of the form
$a_\al\gamma_0$, and denote the integral for such $g$ by $\Psi(\al,f)$.
\np
Let $H$ be the group $PGL(2)$ over $F$. By the Bruhat
decomposition it is $N'A'\cup N'wN'A'$. Define a character on the upper
unipotent subgroup $N'$ by $\psi_{N'}(n')=\psi(x)$, where $n'=n'(x)\in N'$.
\np
For any spherical function $f'$ on $H$ (i.e. $f\in C_c(K'\bksl H/K')$),
define the Fourier orbital integral
$$
\Psi'(\al,f')=\int_{N'}\int_{A'}f'(n'wn'(\al)a')\psi_{N'}(n')\iota(a')da'dn',
$$
where $\iota$ is $\chi_0$ (the unramified quadratic character on $F^\times$)
in the split case, and $1$ in the non-split case.
\np
Let $\pi_\zeta=I_G(\zeta,1/2+\zeta)$ be a (certain) unramified representation of
$G$, induced from its Borel subgroup $B$. Let $\pi_\zeta'=I_H(\zeta,-\zeta)$ be
a (certain) unramified representation of $H$, induced from its Borel subgroup $B'$.
We say that two spherical functions $f$ on $G$ and $f'$ on $H$ are corresponding if
their Satake transforms are equal, i.e.
$\mbox{tr}\ \pi_\zeta(f)=\mbox{tr}\ \pi_\zeta'(f')$,
for all complex numbers $\zeta$.
\np
Consider a pair of corresponding functions $(f,f')$.
The main result of  this paper shows that their
Fourier orbital integrals are related by
$$
\Psi(\al,f)=(\th,\al)\psi(\al)|\al|\Psi'(\al^{-1},f'),
$$
where $(\th,\al)$ is a Hilbert symbol.
\np
Our paper was motivated by Flicker-Mars [FM], which deals with lifting
representations from ${\bf H}=PGL(2)$ to ${\bf G}=PGSp(4)$. It uses the
property that the lifts have periods with respect to the cycle
(closed subgroup) ${\bf C}_\theta$.
Here ${\bf C}_\theta$ is the centralizer of $\dia(a_\theta,a_\theta)$
($a_\theta=\mbox{antidiag}(1,\theta)$) in ${\bf G}$.
This property led [FM] to apply the theory of the Fourier summation formula
for $PGSp(4)$ and the cycle ${\bf C}_\theta$
over a global field. This summation formula is a special case of Jacquet's
relative trace formula. Other approaches to establishing this
lifting of representations use the theory of the Weil representation
(see Oda [O], Rallis-Schiffmann [R], [RS], Langlands [L], Piatetski-Shapiro [PS]).
\np
The Fourier summation formula is obtained by integrating the kernel $K_f(n,h)$,
in fact its product $K_f(n,h)\ovr{\psi}_N(n)$ with the complex conjugate
of the value of the character $\psi_N$ on
${\bf N}(\mathbb{A}$), over $n\in{\bf N}(F)\bksl{\bf N}(\mathbb{A})$ and
$h\in{\bf C}_\theta(F)\bksl{\bf C}_\theta(\mathbb{A})$.
This kernel of the convolution operator $r(f)$ on the space of cusp
forms on ${\bf G}(\mathbb{A})$,
has the geometric expansion $\sum_{\gamma\in {\bf G}(F)}f(n^{-1}\gamma h)$,
and a spectral expansion.
One compares the geometric side of the summation formula on
${\bf G}(\mathbb{A})$ ($=PGSp(4,\mathbb{A})$) to the geometric side of a similar
summation formula on ${\bf H}(\mathbb{A})$ ($=PGL(2,\mathbb{A})$).
The equality of the geometric sides (for different
test functions) implies the equality of the spectral sides
of these two formulae.
This can be used to obtain various conclusions about lifting of representations
from $PGL(2,\mathbb{A})$ to $PGSp(4,\mathbb{A})$.
To carry out the separation argument which plays a key role in these
studies one needs a ``fundamental lemma'', which asserts that corresponding
spherical functions on $PGSp(4)$ and $PGL(2)$ have matching
Fourier orbital integrals.
\np
Proposition 8 of [FM] states a precise form of the conjectured
fundamental lemma. A direct proof of this statement for the unit
elements in the Hecke algebras is given in Proposition 6 of [FM].
\np
Note that under the isomorphism of $PGSp(4)$ with $SO(3,2)$, if
$\th$ is a square in $F^\times$, the image of ${\bf C}_\th$ is
the split group $SO(2,2)$, and if $\theta$ is non-square, it is
$SO(3,1)$. This paper proves the fundamental lemma conjectured in [FM]
for the pairs of local groups $SO(3,2)/SO(3,1)$ and $SO(3,2)/SO(2,2)$;
remarkable cancellations simplify the proof, indeed make it possible.
\np
The proof is based on computing the Fourier transforms of the orbital
integrals (referred to also as the
Mellin transform, since the variable of integration is multiplicative).
By the Fourier inversion formula, the equality of the Fourier
transforms of the orbital integrals implies the equality of the
orbital integrals themselves. This approach avoids dealing with
the asymptotic behaviour of our orbital integrals.
It was first used in Jacquet [J] for the unit element, and then extended
by Mao [M] for the general elements (in their case of $GL(3)$).
In our case, the unit element is treated in [FM] by direct computations.
Here, we give the complete proof of the ``fundamental lemma'' of [FM].
Another interesting question
in this direction is the generalization of these results to the case of
$SO(n)/SO(n-1)$. We hope that this method would apply in this case too.
\np
Our ``Fourier'' situation is significantly different from that of
standard conjugacy, where it is known that: (1) the fundamental lemma for the
unit element implies it for general spherical functions, and (2) the knowledge
of the transfer of orbital integrals of general functions implies in
principle the fundamental lemma for the unit element, and vice versa.
No such results are
known in the ``Fourier'' setting of this paper, in particular since there is
no analogue of the reduction of orbital integrals to ones of elliptic elements
on Levi factors. Of course it will be of much interest to establish analogous
of (1) and (2) in the ``Fourier'' case.
\np
Professor Y. Flicker suggested this problem to me. I
would like to thank him for numerous helpful discussions on this
subject, and Dr. Z. Mao for discussions about his paper.

\noindent
{\bf 0. Statement of results.} Let $F$ be a local non-archimedean
field, of residual characteristic
$\neq 2$. Denote by $R$ the (local) ring of integers of $F$. Let $\pii\in F$
be a generator of the maximal ideal of $R$. Denote by $q$ the number
of elements of the residue field ${\mathbb{F}}=R/\pii R$ of $R$.
Normalize the absolute value on $F^\times$ by $|\pii|=q^{-1}$.
Fix an additive character $\psi$ on $F$ with conductor $R$ (i.e.
$\psi$ is trivial on $R$ but not on $\pii^{-1}R$).
\np
Let $G=SO(3,2;F)$ be the the group of $g\in SL(5;F)$ with $^tgJ g=J$,
where ${}^tg$ is the transpose of $g$ and $J=J_5$. Here $J_n=(\delta_{n+1-i,i})$
is the $n$ by $n$ matrix with 1's on the antidiagonal and $0$'s everywhere
else. Then $G$ is the split special orthogonal group in five variables.
Denote by $V$ the five dimensional vector space
of columns, over $F$. The group $G$ acts on $V$ via multiplication
on the left. In the split case, let $v_0\in V$ be the column $^t(0,0,1,0,0)$.
Set $C=\,\mbox{Stab}_G(v_0)$. Then $C$ is the split special orthogonal
group over $F$ in $4$ variables; we denote it by $SO(2,2;F)$. The symmetric
space $G/C$ is known to be isomorphic
(via the map $g\mapsto gv_0$) to a four dimensional closed
subvariety $S$ of $V$, given by a quadratic equation.
In the non-split case, let
$v_0\in V$ be the column $^t(0,2\th,0,1,0)$, $\th$ not in $(F^\times)^2$.
In Section I.2 we define a subgroup $C'_\th$ of $PGSp(4)$. We denote by
$C_\th=SO(3,1;F)$ the image of $C'_\th$ under the isomorphism
from $PGSp(4)$ to $SO(3,2)$. In Lemma 1.4 we show that
$C_\th=\,\mbox{Stab}_G(v_0)$. The quotient $G/C_\th$ is known to be isomorphic
(via the map $g\mapsto gv_0$) to a four dimensional closed
subvariety $S$ of $V$, given by a quadratic equation.
To simplify the notations, we write $C_\th$ for both split and
non-split cases. The split case corresponds to $C_1$, where $C_1=C$.
\newline\indent
Denote by $P$ the maximal upper triangular parabolic subgroup of $G$
with abelian unipotent radical, $N$.
The subgroup $B$ denotes the upper triangular Borel subgroup of $G$.
Let $K=SO(3,2;R)$ be the standard maximal compact subgroup of $G$. Define
$$
A=\{a_\al=\dia(\al,1,1,1,\al^{-1});\ \al\in F^\times\},\
M=\{\dia(1,m,1);\ m\in SO(J_3)\}.
$$
Then $P=NMA$.
Let $A_0$ be the diagonal subgroup of $G$ and $N_0$ the maximal
unipotent subgroup of $B$.
We have $B=A_0N_0$.
Note that $N$ is a subgroup of $N_0$, isomorphic to $F^3$.
We will write $n=n(x_1,x_2,x_3)$ (as in I.1).
In the split case, define a character $\psi_N$ on $N$
by $\psi_N(n)=\psi(x_2)$, where $n=n(x_1,x_2,x_3)$.
In the non-split case, set $\psi_{N,\th}(n)=\psi(x_1+2\th x_3)$.
\newline\indent
The subgroup $N$ acts on $G/C_\th$ by multiplication on the left,
turning it into a disjoint union of $N$-orbits.
By Proposition 1.5 of Section I the $N$-$C_\th$-double cosets in $G$
of maximal dimension (which is equal to $9$, as $\mbox{dim}(N)=3$ and
$\mbox{dim}(C_\th)=6$)
are represented by $a_\al\gamma_0$, $a_\al\in A$,
for a certain matrix $\gamma_0$, defined at the beginning of Section I.
Moreover the $N$-orbits of $S$ of maximal dimension, $3$,
are of the form $Na_\al\gamma_0v_0$,
$a_\al\in A$.
\newline\indent
As usual, denote by $C(X)$ the space of complex valued functions
on an $l$-space $X$ (see [BZ]), and in $C_c^\infty(X)$, the subscript ``c''
indicates ``compactly supported'', and ``$\infty$'' means ``locally
constant''.
For any $f\in C_c^\infty(G)$, define
$$
\phi_f(gv_0)=\int_C f(gh)dh.
$$
Then $\phi_f\in C_c^\infty(S)$.
If $f$ is a spherical function (i.e. $K$-biinvariant, or
$f\in C_c(K\bksl G/K)$), or even if only $f\in C_c(K\bksl G)$,
then $\phi_f\in C_c(K\bksl S)$.
For any $\phi\in C_c^\infty(S)$, define the orbital integral
$$
\Psi(\al,\phi)=\int_N \phi(na_\al\gamma_0v_0)\psi_N(n)dn.
$$
\newline\indent
Let $H$ be the group $PGL(2,F)$. Its elements will be denoted
by their representatives in $GL(2)$.
Note that $H$ has a trivial center.
Denote by $B'$ the upper triangular
Borel subgroup of $H$. Let $K'=PGL(2,R)$ be the standard maximal
compact subgroup of $H$. We have $B'=N'A'$, where
$$
N'=\left\{n(x)=\Two{1}{x}{0}{1};\ x\in F\right\},\
A'=\left\{\Two{\al}{0}{0}{1};\ \al\in F^\times\right\}.
$$
Define a character $\psi_{N'}$ of $N'$ by
$\psi_{N'}(n(x))=\psi(x)$.
Let $\chi_0$ be the trivial or an unramified quadratic character of $F^\times$
(i.e. $\chi_0(\pii)^2=1$, and $\chi_0(R^\times)=1$).
Put $\iota=\chi_0$ in the split case and $1$ in the non-split case.
Define the integral
$$
\Psi_{H}(\al,f')=\int_{N'}\int_{A'}f'\left(n\Two{0}{1}{-1}{0}
\Two{1}{\al}{0}{1}\Two{a}{0}{0}{1}\right)\iota(a)\psi_{N'}(n)dnd^\times a.
$$
\np
Denote by $C_c(N'\bksl H,\psi_{N'})$ the space of complex valued
compactly supported modulo $N'$ functions $\phi'$ on $H$,
which satisfy (for any $n\in N'$) the relation
$$
\phi'(ng)=\ovr{\psi}_{N'}(n)\phi'(g),
$$
where $\ovr{z}$ denotes the complex conjugate of $z$.
Write $\Ch$ for the space of such right $K'$-invariant functions.
Given $f'\in C_c^\infty(H)$, define a function
$\phi_{f'}'\in C_c(N'\bksl H,\psi_{N'})$ on $H$ by
$$
\phi_{f'}'(g)=\int_{N'}\psi_{N'}(n)f'(ng)dn.
$$
If $f'$ is $K'$-biinvariant, then $\phi'_{f'}\in\Ch$.
Define the integral
$$
\Psi'(\al,\phi_{f'}')=\int_{A'}\phi'_{f'}
\left(\Two{0}{1}{-1}{0}\Two{1}{\al}{0}{1}\Two{a}{0}{0}{1}\right)
\iota(a)d^\times a.
$$
Thus, by definition $\Psi'(\al,\phi_{f'}')=\Psi_{H}(\al,f')$.

{\bf Definition.} The functions $f\in C_c^\infty(G)$ and
$f'\in C_c^\infty(H)$ are called {\it matching} if for every $\al\in F^\times$
we have
$$
\Psi(\al,\phi_f)= (\th,\al)\psi(\al)|\al|\Psi'(\al^{-1},\phi_{f'}').
$$
Note that in the split case $(\th,\al)=1$.
\np
Let $\pi=I_G(\zeta,\zeta')$ be the representation of the group $G$ which is
normalizedly induced from the character $|\al_1|^\zeta|\al_2|^{\zeta'}$
of the Borel subgroup, $B$, where $\al_1$ and $\al_2$ are the two simple
roots of $G$ with respect to $B$.
The space of this representation consists of locally constant functions $\phi$,
such that
$$
\phi(nak)=\delta_B^{1/2}(a)|\al_1(a)|^\zeta|\al_2(a)|^{\zeta'}\phi(k),
$$
where $n\in N$, $a=\dia(\al,\beta,1,\beta^{-1},\al^{-1})$ and
$\al_1(a)=\al/\beta$, $\al_2(a)=\beta$; $G$ acts by right translation.
Let $f$ be a $K$-biinvariant,
compactly supported function. Define its Satake transform $f^\vee$ by
$f^\vee(\pi)=\,\mbox{tr}\,\pi(f)$.
\np
Let $\pi_\zeta'=I_{H,\chi_0}(\zeta,-\zeta)$ be the
representation of the group $H$
which is normalizedly induced from the character
$$
\Two{a}{n}{0}{b}\longmapsto \left|\frac{a}{b}\right|^{\zeta}\chi_0
\biggl(\frac{a}{b}\biggr)
$$
of $B'$.
\np
Let $f'$ be a $K'$-biinvariant, compactly
supported function on $H$. Its Satake transform is defined again
by $f'^\vee(\pi')=\,\mbox{tr}\,\pi'(f')$.

{\bf Definition.} The $K$-biinvariant function $f$ on $G$ and the
$K'$-biinvariant function $f'$ on $H$ are called {\it corresponding} if
for any complex number $\zeta$ we have
$$
f^\vee(\pi_\zeta)=f'^\vee(\pi_\zeta'),
$$
where $\pi_\zeta=I_G(\zeta,1/2+\zeta)$ and
$\pi_\zeta'=I_{H,\chi_0}(\zeta, -\zeta)$ are the representations of $G$
and $H$ defined above.

The unit elements $f^0$ and $f'^0$ of the Hecke algebras $C_c(K\bksl G/K)$
and $C_c(K'\bksl H/K')$, which are the characteristic functions of $K$
and $K'$ divided by their volumes, are corresponding.
It is shown in [FM] that they are
matching. The main result of this paper is the following extension of
that result, conjectured in [FM].

\noindent
{\bf Theorem.} {\em Corresponding $f$ and $f'$ are matching. }

Our approach is analogous to that of [J], [M].
An alternative approach would be to directly compute the integral.
The general structure of the proof is as follows.
Under the action of $K$, $S$ can be decomposed into $K$-orbits.
Each orbit has a representative (see Proposition 1.7)
of the form $d_rv_1$ ($r\geq 0$), where
$d_r=\dia(\pii^r,1,1,1,\pii^{-r})$. In the split case
$v_1={}^t(1/2,0,0,0,1)$, in the non-split case $v_1={}^t(2\th,0,0,1,1)$.
The main result of Section I is Proposition 1.8, which computes
the volume of $Kd_rv_1$. This section contains also some results needed
in Section II.
\np
We write ${\cal F}(\phi)=\phi'$
(where $\phi\in C_c(K\bksl S)$ and $\phi'\in\Ch$) if
$$
\int_S\phi(s)T_\zeta(s)ds=\int_{A'}\phi'(a)W_\zeta(a)|a|^{-1}da.
$$
Here $T_\zeta$ is the $K$-invariant function on $S$, such that
$T_\zeta(d_0v_1)=1$ and for $r\geq 1$ (see Proposition 2.2)
$$
T_\zeta(d_rv_1)=\sum_{\xi\in\{\zeta,-\zeta\}}
\frac{(q^{\frac{3}{2}+\xi}-1)(1\mp q^{-\frac{1}{2}-\xi})}
{(q^{\xi}-q^{-\xi})(q^{\frac{3}{2}}\mp q^{-\frac{1}{2}})}
q^{-r(\frac{3}{2}-\xi)},
$$
where the ``$-$'' sign occurs in the split case and ``$+$''
in the non-split case.
The function $W_\zeta$ is the normalized unramified Whittaker function in
the space of the representation $\pi_\zeta'$ (see Proposition 2.4).
\np
We show in Proposition 2.5 that ${\cal F}$ is a linear bijection
between the spaces $C_c(K\bksl S)$ and $\Ch$.
\np
The map $f\mapsto \phi_f$ from $C_c(K\bksl G/K)$ to $C_c(K\bksl S)$,
and the map $f'\mapsto \phi'_{f'}$ from
$C_c(K'\bksl H/K')$ to $\Ch$,
are used in Section II to show that the relation
$f^\vee(\pi_\zeta)=f'^\vee(\pi_\zeta')$
is equivalent to ${\cal F}(\phi_f)=\phi'_{f'}$.
\np
Define $\phi_r$ to be the characteristic function of the orbit $Kd_rv_1$.
Since $K\bksl S$ is the disjoint union of $Kd_rv_1$, $r\geq 1$, (Proposition 1.3)
the set $\{\phi_r;r\geq 0\}$ is a basis of the space $C_c(K\bksl S)$.
Define $\Phi_r=\sum_{i=0}^r\phi_i$. Then the set $\{\Phi_r;r\geq 0\}$ is also
a basis of $C_c(K\bksl S)$.
\np
Now, we consider the group $H$.
Since $H=N'A'K'$, any function in $C_c(N'\bksl H/K',\psi_{N'})$ is
defined by its values on $A'$. For $r\geq 0$, define the function $\phi_r'$
in this space by
$$
\phi_r'\left(\Two{\al}{0}{0}{1}\right)=\left\{\begin{array}{cl} 1, &
\ \mbox{ if } |\al|=q^{-r}, \\ 0,&\ \mbox{ otherwise }. \end{array}\right.
$$
We show in Proposition 2.4(2) that $\phi'(\dia(\al,1))=0$ if $|\al|>1$.
Hence, the set $\{\phi_r';r\geq 0\}$ is a basis of $\Ch$.
Without lost of generality we assume that $\chi_0(\pii)=-1$. Indeed, if this
is not the case then we can change the basis $\phi_r'\mapsto (-1)^r\phi_r'$.
Set $\phi_r'=0$ if $r<0$. The main result of Section II asserts
that for any integer $r\geq 0$, we have ${\cal F}(\Phi_r)=
(-1)^rq^r(\phi_r'\pm\phi_{r-1}')$, i.e.
$$
\int_S\Phi_r(s)T_\zeta(s)ds=(-1)^rq^r\int_{A'}(\phi_r'(a)\pm
\phi_{r-1}'(a))W_\zeta(a)|a|^{-1}da,
$$
where as usual the ``$+$'' sign occurs in the split case, and
the ``$-$'' in the non-split case.
Thus, if two spherical functions $f\in C_c^\infty(G)$ and
$f'\in C_c^\infty(H)$ are corresponding we have
${\cal F}(\phi_f)=\phi'_{f'}$.
Since ${\cal F}$\,:\,$C_c(K\bksl S)\rightarrow\Ch$ is an isomorphism
of two vector spaces, to prove that corresponding functions are matching (i.e.
$\Psi(\al,\phi_f)=(\th,\al)\psi(\al)|\al|\Psi'(\al^{-1},\phi'_{f'})$)
it is enough to show that (for $r\geq 0$)
$$
\Psi(\al,\Phi_r)=(\th,\al)\psi(\al)|\al|
\Psi'(\al^{-1},(-1)^rq^r(\phi'_r\pm\phi'_{r-1})).
$$
\np
In Section III we show that (for any $r\geq 1$)
$$
\int_{F^\times}\Psi(\al,\Phi_r)\chi(\al)d^\times\al=
\int_{F^\times}(\th,\al)\psi(\al)|\al|\Psi'(\al^{-1},(-1)^rq^r
(\phi'_r\pm\phi'_{r-1}))\chi(\al)d^\times\al,
$$
where $\chi$ is any complex valued character of $F^\times$.
The case $r=0$ follows from this result and from the case of the unit
element, treated in [FM].
If $\chi$ is ramified both integrals are equal to $0$.
The Fourier inversion formula now implies the required result for
the split case.

\begin{center}
{\bf I. The group $G$, subgroup $C$ and the $K$-orbits of $G/C$.}
\end{center}
{\bf I.1. The group $G$.} The group ${\bf G}=SO(3,2)$ can also be
defined as
$$
\{g\in GL(5);\, Q(gv,gv)=Q(v,v),\, \det(g)=1 \},
$$
where $Q(v,w)={}^tvJw$, hence $Q(v,v)={}^tvJv=2v_1v_5+2v_2v_4+v_3^2$
is a quadratic form on the 5 dimensional vector space $V$ of columns.
Let ${\bf P}$ be the maximal upper triangular parabolic
subgroup of ${\bf G}$ with abelian unipotent radical, ${\bf N}$.
Let ${\bf B}$ be the upper triangular Borel subgroup of ${\bf G}$.

{\bf Definition.} Define the matrix
\EQ\label{1.0}
n=n(x_1,x_2,x_3,x_4)=
\left(\ColF{1}{0}{0}{0}{0}\ColF{x_3'}{1}{0}{0}{0}\ColF{x_2'}
{-x_4}{1}{0}{0}\ColF{x_1'}{-\frac{1}{2}x_4^2}{x_4}{1}{0}
\ColF{z}{x_1}{x_2}{x_3}{1}\right),
\EN
where
\EQ\label{1.1}
x_3'+x_3=0,\ x_2+x_2'=x_3x_4,\ x_1+x_1'=-x_2x_4+\frac{1}{2}x_3x_4^2,\
z=-x_1x_3-\frac{1}{2}x_2^2.
\EN
Let ${\bf A}_0=\{ \dia(\al,\beta,1,\beta^{-1},\al^{-1});\ \al,\beta\neq 0\}$
be the diagonal subgroup of ${\bf G}$. Let ${\bf N}_0=\{ n=n(x_1,x_2,x_3,x_4)\}$
be the upper triangular maximal unipotent subgroup of ${\bf B}_0$.
Put $n(x_1,x_2,x_3)=n(x_1,x_2,x_3,0)$.
We have ${\bf B}={\bf A}_0{\bf N}_0$ and ${\bf P}={\bf NMA}$, where
$$
{\bf M}=\left\{\Three{1}{0}{0}{0}{m}{0}{0}{0}{1};\
m\in SO(J_3)\right\}, \ \
{\bf A}=\left\{\Three{\al}{0}{0}{0}{I}{0}{0}{0}{\al^{-1}};\ \al\neq 0\right\},
$$
and ${\bf N}=\{ n=n(x_1,x_2,x_3) \}$.
The standard Levi subgroup of ${\bf P}$ is the product ${\bf MA}$.
\np
If $x_4=0$ the condition (\ref{1.1}) reduces to
$x_i'+x_i=0$ ($i=1,2,3$) and $z=-x_1x_3-\frac{1}{2}x_2^2$.
We define $G={\bf G}(F)$, $P={\bf P}(F)$, $N={\bf N}(F)$,
$A={\bf A}(F)$, $N_0={\bf N}_0(F)$ and $A_0={\bf A}_0(F)$.
Define the character $\psi_N$ on $N$.
In the split case, let $\psi_N(n(x_1,x_2,x_3))=\psi(x_2)$.
In the non-split case, let $\psi_{N,\th}(n(x_1,x_2,x_3)) =
\psi(x_1+2\th x_3)$.
\np
Consider the split case. Put
$$
{\bf C}=\left\{ \Three{A_1}{0}{A_2}{0}{1}{0}{A_3}{0}{A_4};\,
\Two{A_1}{A_2}{A_3}{A_4}\in SO(J_4)\right\}.
$$

\noindent{\bf Lemma 1.1.} {\em  Put $v_0={}^t(0,0,1,0,0)$, and $C={\bf C}(F)$.
Then $C$ is the stabilizer of $v_0$ under the action of $G$ on $V$, and the map
$g\mapsto gv_0$ embeds $G/C$ into $S$, where
$S$ is the sphere $v$ in $V$ such that $Q(v,v)=1$. }

\noindent{\em Proof.} Clearly $C=\mbox{Stab}_G(v_0)$.
Since $G$ is the group $SO(J)$, and the third
column $x$ of any element $g$ of $G$ is $gv_0$, $x$
satisfies the condition $Q(x,x)=1$.
\hspace*{\fill}$\Box$

\noindent{\bf Remark.} Note that
$$
S=\{{}^t(x_1,x_2,x_3,x_4,x_5)\in V;\ 2x_1x_5+2x_2x_4+x_3^2=1\}.
$$

\noindent
{\bf I.2. The isomorphism between $SO(3,2)$ and $PGSp(4)$.}
Let $G'$ be the group $PGSp(4,F)$ of matrices
$g\in GL(4,F)$ such that $gJ'^tg=\lam J'$, where $J'$ is the matrix
$\ant(1,1,-1,-1)$ and $\lam\in F^\times$.
Fix $\th\in F^\times$ which is not a square.
Let $a_\th=\ant(1,\th)$. Let $C'_\th$ be the
centralizer of $\dia(a_\th,a_\th)$ in $G'$. Let $N'$ be the unipotent radical
of the Siegel parabolic subgroup $P'$ of type $(2,2)$ of $G'$. Recall that
$$
N'=\left\{ n=\Two{I}{A}{0}{I}\in G';\ A=\Two{x}{y}{z}{x}\right\}.
$$
Fix a complex valued non-trivial character $\psi$ of $F$ and define
the character $\psi_\th$ of $N'$ by $\psi_\th(n)=\psi(z-\th y)$.
The stabilizer of this character is a non-split torus (see [FM]).

\indent
{\bf Definition.} Define a five dimensional space $X$ by
$$
X=\{T\in M_4(F);\ {}^t(TJ')=-TJ', \mbox{tr}(T)=0\}.
$$
Choose the basis $\{e_1,e_2,e_3,e_4,e_5\}$ of $X$ so that
$$
T=T(x_1,x_2,x_3,x_4,x_5)=x_1e_1+x_2e_2+x_3e_3+x_4e_4+x_5e_5
$$
is represented by the matrix
\EQ\label{12.1}
\left(\CFour{-x_3/2}{x_2/2}{x_5}{0}\CFour{x_4}{x_3/2}{0}{-x_5}
\CFour{x_1/2}{0}{x_3/2}{x_2/2}\CFour{0}{-x_1/2}{x_4}{-x_3/2}\right).
\EN
The inner form on this space is $(T_1,T_2)=\mbox{tr}(T_1T_2)$, where
$T_1$, $T_2$ are in $X$. Define an action of $G'$ on $X$ via
$g$: $T\mapsto gTg^{-1}$.

\noindent
{\bf Lemma 1.2.} {\em The action $g$: $T\mapsto gTg^{-1}$ of
$G'=PGSp(4)$ on $X$ is well defined and establishes an isomorphism
from $G'=PGSp(4)$ to $G=SO(3,2)=SO(J)$. }

\noindent
{\em Proof.} To show that this action is well defined, we have to
prove that ${}^t(gTg^{-1}J')=-gTg^{-1}J'$.
This relation is equivalent to
${}^tJ'{}^tg^{-1}{}^tT{}^tg=-gTg^{-1}J'$.
Multiplying both sides by $g^{-1}$ on the left and by ${}^tg^{-1}$
on the right, we obtain
$g^{-1}{}^tJ'{}^tg^{-1}{}^tT=-Tg^{-1}J'{}^tg^{-1}$.
But $g^{-1}J'{}^tg^{-1}=\lambda J'$ implies that
$g^{-1}{}^tJ'{}^tg^{-1}=\lambda{}^tJ'$.
We arrive at
$$
\lambda{}^tJ'{}^tT=-Tg^{-1}J'{}^tg^{-1}=-\lambda TJ',
$$
which is true since $T\in X$.
\np
Further if $T_1=T_1(x_1,x_2,x_3,x_4,x_5)$
and $T_2=T_2(y_1,y_2,y_3,y_4,y_5)$ then
$$
\mbox{tr}(T_1T_2)=x_1y_5+x_2y_4+x_3y_3+x_4y_2+x_5y_1.
$$
Since $\mbox{tr}(gT_1T_2g^{-1})=\mbox{tr}(T_1T_2)$, the action of $G'$ on $X$
defines an orthogonal group on $X$. The space $X$ is isomorphic
to the $5$ dimensional vector space $V$, from Section I.1.
Thus this orthogonal group is the group $G=SO(J)$.
\hspace*{\fill}$\Box$

\noindent
{\bf Lemma 1.3.} {\em Under the isomorphism of Lemma 1.1, the image of
subgroup $C'_\th$ is $C_\th$, the centralizer (in $G$) of
$$
\left(\ColF{-1}{0}{0}{0}{0}\ColF{0}{0}{0}{\frac{1}{2}\th^{-1}}{0}
\ColF{0}{0}{-1}{0}{0}\ColF{0}{2\th}{0}{0}{0}\ColF{0}{0}{0}{0}{-1}\right).
$$ }

\noindent
{\em Proof.} By matrix multiplication
$$
\left(\CFour{0}{\th}{0}{0}\CFour{1}{0}{0}{0}\CFour{0}{0}{0}{\th}
\CFour{0}{0}{1}{0}\right)
\left(\CFour{-x_3/2}{x_2/2}{x_5}{0}\CFour{x_4}{x_3/2}{0}{-x_5}
\CFour{x_1/2}{0}{x_3/2}{x_2/2}\CFour{0}{-x_1/2}{x_4}{-x_3/2}\right)
\left(\CFour{0}{1}{0}{0}\CFour{\th^{-1}}{0}{0}{0}\CFour{0}{0}{0}{1}
\CFour{0}{0}{\th^{-1}}{0}\right)
$$
$$
=\left(\CFour{x_3/2}{\th x_4}{-x_5}{0}\CFour{x_2/(2\th)}
{-x_3/2}{0}{x_5}\CFour{-x_1/2}{0}{-x_3/2}
{\th x_4}\CFour{0}{x_1/2}{x_2/(2\th)}{x_3/2}\right),
$$
which implies the lemma.
\hspace*{\fill}$\Box$

Note that under the isomorphism between
$G'$ and $G$, the unipotent subgroup $N'$ of $G'$ is isomorphic to
the subgroup $N$ of $G$ via
$$
\left(\CFour{1}{0}{0}{0}\CFour{0}{1}{0}{0}\CFour{x}{z}{1}{0}\CFour{y}{x}{0}{1}
\right)\mapsto
\left(\ColF{1}{0}{0}{0}{0}\ColF{-y}{1}{0}{0}{0}\ColF{-2x}
{0}{1}{0}{0}\ColF{2z}{0}{0}{1}{0}\ColF{2zy-2x^2}{-2z}{2x}{y}{1}\right).
$$
In particular $z-\th y\mapsto -\frac{1}{2}(x_1+2\th x_3)$ which justifies
the choice of the character $\psi_{N,\th}$ on $N$.

\noindent
{\bf Lemma 1.4.} {\em Put $v_0={}^t(0,2\th,0,1,0)$. Then $C_\th$ is the
stabilizer of $v_0$ under the action of $G$ on $V$, and the map
$g\mapsto gv_0$ embeds $G/C_\th$ into $S$, where $S$ is the sphere $v$ in $V$
such that $Q(v,v)=4\th$. }

\noindent
{\em Proof.} The image of the subgroup $C'_\th$ in $G$ is the subgroup $C_\th$,
which consists of matrices of the form
$$
\left(\ColF{a_{11}}{a_{21}}{a_{31}}{-a_{21}/(2\th)}{a_{51}}
\ColF{a_{12}}{a_{22}}{a_{32}}{(1-a_{22})/(2\th)}{a_{52}}
\ColF{a_{13}}{a_{23}}{a_{33}}{-a_{23}/(2\th)}{a_{53}}
\ColF{-2\th a_{12}}{2\th(1-a_{22})}{-2\th a_{32}}{a_{22}}{-2\th a_{52}}
\ColF{a_{15}}{a_{25}}{a_{35}}{-a_{25}/(2\th)}{a_{55}}\right).
$$
Clearly $C_\th=\mbox{Stab}_G(v_0)$. Recall that $Q(v,w)={}^tvJw$.
If ${\bf y}_2$ is the second and ${\bf y}_4$
is the fourth columns of the orthogonal group $G$, then they satisfy $
Q({\bf y}_2,{\bf y}_2)=Q({\bf y}_4,{\bf y}_4)=0$ and
$Q({\bf y}_2,{\bf y}_4)=Q({\bf y}_4,{\bf y}_2)=1$.
The element $gv_0$ is the sum $2\th{\bf y}_2+{\bf y}_4$. Hence, we have
$$
Q(2\th{\bf y}_2+{\bf y}_4,2\th{\bf y}_2+{\bf y}_4)=4\th^2 Q({\bf y}_2,{\bf y}_2)
+4\th Q({\bf y}_2,{\bf y}_4)+Q({\bf y}_4,{\bf y}_4)=4\th.
$$
\hspace*{\fill}$\Box$

{\bf Remark.} The sphere $S$ is equal to
$$
S=\{{}^t(x_1,x_2,x_3,x_4,x_5)\in V;\ 2x_1x_5+2x_2x_4+x_3^2=4\th\}.
$$

\noindent
{\bf I.3. Double coset decomposition.}
In the split case, we define
$$
v_1=\ColFive{\frac{1}{2}}{0}{0}{0}{1},\
\gamma_0=\left(\ColF{\frac{1}{2}}{0}{-1}{0}{-1}\ColF{0}{1}{0}{0}{0}
\ColF{\frac{1}{2}}{0}{0}{0}{1}\ColF{0}{0}{0}{1}{0}
\ColF{-\frac{1}{4}}{0}{-\frac{1}{2}}{0}{\frac{1}{2}}\right),\
v_0=\ColFive{0}{0}{1}{0}{0}.
$$
In the non-split case ($\th\not\in (F^\times)^2$), we define
$$
v_1=\ColFive{2\th}{0}{0}{1}{1},\
\gamma_0=\left(\ColF{1}{-1}{0}{0}{0}\ColF{1}{0}{0}{0}{0}
\ColF{0}{0}{1}{0}{0}\ColF{0}{0}{0}{1}{1}\ColF{0}{0}{0}{-1}{0}\right),\
v_0=\ColFive{0}{2\th}{0}{1}{0}.
$$
Note that $v_1=\gamma_0v_0$.

\noindent{\bf Proposition 1.5.} {\em We have the disjoint decomposition:
$G=PC_\th\cup NA\gamma_0C_\th$.
The representatives of the $N$-$C_\th$-orbits of maximal dimension,
which is $9$, are of the form $a\gamma_0$, $a\in A$. Futhermore,
the map $g\mapsto gv_0$ of Lemmas 1.1 and 1.3 establishes an isomorphism
of homogeneous spaces from $G/C_\th$ to $S$.}

\noindent{\em Proof.}
Consider the left action of $P$ on $S$. Since the last row of any
element of $P$ is of the form $(0,0,0,0,\alpha^{-1})$, $\alpha\neq 0$,
we conclude that there are at least two $P$-invariant subsets in $S$:
a closed subset
$\{ x={}^t(x_1,x_2,x_3,x_4,0);\ Q(x,x)=1\}$, and an open subset
$\{ x={}^t(x_1,x_2,x_3,x_4,x_5);\ Q(x,x)=1,\ x_5\neq 0\}$.
We claim that $P$ acts transitively on each of these two subsets.
\np
Consider the split case. The element $v_1={}^t(1/2,0,0,0,1)$, which is a
representative of the open subset. Acting first by an element from $A$
and then from $N$, we obtain the transpose of
$$
\left(\frac{\al}{2}+\frac{z}{\al},\frac{x_1}{\al},\frac{x_2}{\al},
\frac{x_3}{\al},\frac{1}{\al}\right).
$$
When $(x_1,x_2,x_3)$ runs through $F^3$ and $\alpha$ over $F^\times$, this
column runs through all elements ${}^t(x_1,x_2,x_3,x_4,x_5)$
of $S$, with $x_5\neq 0$. Thus $P$ acts transitively on the open $P$-subset
of $S$, i.e. it is a $P$-orbit.
\np
Consider the non-split case. The element $v_1={}^t(2\th,0,0,1,1)$
is a representative of the open subset. Acting first by an element from $A$
and then from $N$, we obtain the transpose of
$$
\left(2\al\th-x_1+\frac{z}{\al},\frac{x_1}{\al},
\frac{x_2}{\al},1+\frac{x_3}{\al},\frac{1}{\al}\right).
$$
When $(x_1,x_2,x_3)$ runs through $F^3$ and $\alpha$ over $F^\times$, this
column runs through all elements ${}^t(x_1,x_2,x_3,x_4,x_5)$
of $S$, with $x_5\neq 0$. Thus $P$ acts transitively on the open $P$-invariant
subset of $S$, making it a $P$-orbit.
\np
Consider $v_0={}^t(0,0,1,0,0)$ in the split case and $v_0={}^t(0,2\th,0,1,0)$
in the non-split case. This is a representative of the closed subset.
First, acting by an element $n(0,-s,0)\in N$ (see I.1), in the split case,
and $n(-s,0,0)$ in the non-split case,
we can send $v_0$ to ${}^t(s,0,1,0,0)$
(respectively ${}^t(s,2\th,0,1,0)$), $s$ arbitrary. Multiplying it
on the left by an element $g=\dia(1,m,1)\in M$,
we obtain ${}^t(s,m_{21},m_{22},m_{23},0)$,
where ${}^t(m_{21},m_{22},m_{32})$ is the second column of $m$.
Thus, we obtain all elements ${}^t(x_1,x_2,x_3,x_4,0)$ $\in S$.
Note that the $N$-orbits in $S$ of such elements are of dimension $1$.
\np
Consequently, the decomposition of $G$ with respect to
$P$ and $C_\th$ is $G=PC_\th\cup P\gamma_0 C_\th$. Recall that $P=NAM$.
We assert that $\gamma_0^{-1}M\gamma_0\subset C_\th$. Indeed,
from $M\subset\mbox{Stab}(v_1)$ it follows that
$\gamma_0^{-1}M\gamma_0\subset\mbox{Stab}(v_0)$, since $v_1=\gamma_0v_0$.
Thus, $G=PC\cup NA\gamma_0C$.
\np
Since $A\subset C_\th$, we have $PC_\th=NMC_\th$. We have seen that
the $N$-orbits of $NMC_\th/C_\th$ are of dimension $1$. Hence
the $N$-$C_\th$-double cosets of $PC_\th$ are of dimension $7$.
\np
We have seen above that
the map $g\mapsto gv_0$ is an onto map with kernel $C_\th$.
Consequently $G/C_\th\simeq S$. The proposition follows.
\hspace*{\fill}$\Box$

\noindent
{\bf I.4. The subset $Kb_1K$.} Set $b_1=\dia(\pii^{-1},1,1,1,\pii)$,
$b_2=\dia(1,\pii^{-1}\!,1,\pii,1)$ and consider the double coset $Kb_1K$.
Recall that $\mathbb{F}$ is the residue field of $F$,
i.e. a finite field of $q$ elements, $q$ is odd. Define
${\mathbb{N}}_0={\bf N_0}(\mathbb{F})$.
More explicitly
$$
{\mathbb{N}}_0=\left\{n(x_1,x_2,x_3,x_4)=
\left(\ColF{1}{0}{0}{0}{0}\ColF{x_3'}{1}{0}{0}{0}\ColF{x_2'}
{-x_4}{1}{0}{0}\ColF{x_1'}{-\frac{1}{2}x_4^2}{x_4}{1}{0}
\ColF{z}{x_1}{x_2}{x_3}{1}\right);\ x_1,x_2,x_3,x_4\in \mathbb{F}\right\},
$$
where $z$, $x_4$, $x_3'$, $x_3$, $x_2'$, $x_2$, $x_1'$, $x_1$ satisfy the
relations (\ref{1.1}). Define in ${\mathbb{N}}_0$ the subgroups
(the order of $N_i$ is $q^{4-i}$):
$$
\begin{array}{ccc}
N_1 & = &\{n\in{\mathbb{N}}_0;\,n=n(x_1,x_2,x_3,0)\},     \\
N_2 & = &\{n\in{\mathbb{N}}_0;\,n=n(x_1,0,0,x_4)\},       \\
N_3 & = &\{n\in{\mathbb{N}}_0;\,n=n(0,0,x_3,0)\}.
\end{array}
$$
We regard $N_i$ and ${\mathbb{N}}_0$ as subsets of ${\bf N}_0(\mathbb{R})$
on choosing representatives in $R$ of the elements of ${\mathbb{F}}=R/\pii$.

The following Proposition is used in the proof of Proposition 1.8.

\noindent
{\bf Proposition 1.6.} {\em We have the disjoint decomposition
\EQ\label{l1.2}
Kb_1K=Kb_1N_1\cup Kb_2N_2\cup Kb_2^{-1}N_3\cup Kb_1^{-1}.
\EN }
\newline\noindent
{\em Proof.} The Weyl group
$W=\{1,(15),(24),(12)(45),(14)(25),(15)(24),(1452),(1254)\}$
embeds in $K$, so we have
$$
Kb_1K=Kb_1^{-1}K=Kb_2K=Kb_2^{-1}K.
$$
Let $W_1$ be the subset $\{1,(15),(12)(45),(14)(25)\}$ of $W$, and
$P_1=K\cap b_1Kb_1^{-1}$ a parahoric (see, e.g., [T])
subgroup of $K$. Using the Iwahori decomposition:
$$
K=P_1W_1P_1=\cup_{w\in W_1}P_1wP_1,
$$
we have
$$
Kb_1^{-1}K=\cup_{w\in W_1}Kb_1^{-1}P_1wP_1.
$$
Since $b_1^{-1}P_1b_1\in K$, this is equal to
$$
\cup_{w\in W_1}Kb_1^{-1}wP_1.
$$
Since $W_1\subset K$, and $\{w^{-1}b_1^{-1}w; w\in W_1\}$ is the set
$\{b_1^{-1},b_1,b_2^{-1},b_2\}$, we obtain
$$
Kb_1^{-1}P_1\cup Kb_1P_1\cup Kb_2^{-1}P_1\cup Kb_2P_1.
$$
\indent
In our analysis below, we use the decomposition $P_1=N'M_KA_KN_K$,
where $N'$ is the subgroup of ${}^tN$ with underdiagonal entries
from $\pii R$, $A_K=A\cap K$, $N_K=N\cap K$ and
$M_K=M\cap K$ is the maximal compact of $M=\{\dia(1,m,1); m\in SO(J_3)\}$.
\np
To describe the four double cosets, we introduce:
\newline\indent
{\em Case of $Kb_1^{-1}P_1$.}
Since $b_1^{-1}P_1b_1\subset K$, we have $Kb_1^{-1}P_1=Kb_1^{-1}$.
\newline\indent
{\em Case of $Kb_1P_1$.} Since $b_1N'b_1^{-1}\subset K$ and
$M_KA_K$ commutes with $b_1$, we have $Kb_1P_1=Kb_1N_K$.
\np
The set of elements $n\in N_K$ with entries above the diagonal
from $\pii R$ is a normal subgroup
of $N_K$ of elements satisfying $b_1nb_1^{-1}\in K$.
Since the quotient of $N_K$ by this subgroup is $N_1$, we have
$Kb_1N_K=Kb_1N_1$.
\newline\indent
{\em Case of $Kb_2^{-1}P_1$.}
Since $b_2^{-1}N'b_2\subset K$, we have $Kb_2^{-1}P_1=Kb_2^{-1}M_KN_K$.
\np Let $W_2$ be the subgroup of two elements $\{1,(13)\}$,
where $(13)$ is represented by the matrix $w_2=\ant(1,-1,1)$, and put
$P_2=M_K\cap b_2M_Kb_2^{-1}$.
Using the Iwahori decomposition $M_K=P_2W_2P_2$, our set is
$Kb_2^{-1}P_2W_2P_2N_K$. Since $b_2^{-1}P_2b_2\in K$, this is
$Kb_2^{-1}W_2P_2N_K$.
But $W_2=\{1,w_2\}$ and $w_2b_2^{-1}w_2=b_2$, so this double
coset is
$$
Kb_2^{-1}P_2N_K\cup Kb_2^{-1}w_2P_2N_K=Kb_2^{-1}N_K\cup Kb_2P_2N_K.
$$
\indent
{\em Subcase of $Kb_2^{-1}N_K$.}
If $n\in N_K$ is such that $b_2^{-1}nb_2\in K$, then
$n=n(x_1,x_2,x_3,x_4)$, where $x_1,x_2,x_4\in R$ and $x_3\in \pii R$.
The set of such elements is a normal subgroup of $N_K$.
Since the quotient of $N_K$ by this subgroup is $N_3$, we obtain
$$
Kb_2^{-1}N_K=Kb_2^{-1}N_3.
$$
\np{\em Subcase of $Kb_2P_2N_K$.}
To simplify the notations, write $m\in GL(3)$ for $\dia(1,m,1)\in GL(5)$.
Then $P_2=U'A_2U$, where we put
$A_2=\{\dia(\al,1,\al^{-1}); \al\in F^\times, |\al|=1\}$,
$$
U'=\left\{\Three{1}{0}{0}{-x}{1}{0}{-\frac{1}{2}x^2}{x}{1}\!;\ x\in
\pii R\right\},
\ U=\left\{\Three{1}{-y}{-\frac{1}{2}y^2}{0}{1}{y}{0}{0}{1}\!;\ y\in R\right\}.
$$
Since $b_2U'b_2^{-1}\subset K$
and $A_2$ commutes with $b_2$, the coset $Kb_2P_2N_K$ is equal to
$$
Kb_2U'A_2UN_K=Kb_2UN_K.
$$
Note that $UN_K=N_0\cap K$. Finally,
any $n\in N_0$ such that $b_2nb_2^{-1}\in K$, is of the form
$n=n(x_1,x_2,x_3,x_4)$, where $x_2,x_3\in R$ and $x_1,x_4\in \pii R$.
The set of such elements is a subgroup of ${\bf N}_0(R)$, and
$N_2$ is the set of representatives of its left cosets in
${\bf N}_0(R)$. Thus $Kb_2P_2N_K=Kb_2N_2$.
\np{\em Case of $Kb_2P_1$.} Decomposing $P_1$ and using
$b_2N'b_2^{-1}\subset K$ and $b_2A_Kb_2^{-1}=A_K\subset K$, we obtain
$$
Kb_2P_1=Kb_2N'A_KM_KN_K=Kb_2M_KN_K.
$$
Applying the Iwahori decomposition to $M_K$, ($W_2=\{1,w_2\}$), this is
$$
Kb_2P_2W_2P_2N_K=Kb_2P_2N_K\cup Kb_2P_2w_2P_2N_K.
$$
The double coset $Kb_2P_2N_K$ has already been considered.
We are left with the double coset $Kb_2P_2w_2P_2N_K$.
\np
Since $P_2=U'A_2U$ and $w_2U'A_2w_2\subset P_2$, we have
$$
P_2w_2P_2=P_2w_2U'A_2U=P_2w_2U.
$$
Futhermore, since $b_2U'A_2b_2^{-1}\subset K$, the double coset
$Kb_2P_2w_2P_2N_K$ is equal to
$$
Kb_2P_2w_2UN_K=Kb_2U'A_2Uw_2UN_K=Kb_2Uw_2UN_K.
$$
Since $w_2^{-1}b_2w_2=b_2^{-1}$ this is equal to
$$
Kb_2w_2w_2^{-1}Uw_2UN_K=Kb_2^{-1}U_1UN_K,\ \mbox{ where }
U_1=\left\{\Three{1}{0}{0}{-y}{1}{0}{-\frac{1}{2}y^2}{y}{1};\ y\in \Bbb{F}
\right\}.
$$
According to the Iwasawa decomposition, any element $n\in U_1$
can be written as $n=u_1sw_2u_2$, where $u_i\in U$, and $s$, the diagonal
matrix, can be commuted across $b_2^{-1}$. Since $w_2^{-1}b_2^{-1}w_2=b_2$,
we have
$$
Kb_2^{-1}U_1UN_K=Kb_2UN_K,
$$
obtaining again a coset which has already been considered.
\hspace*{\fill}$\Box$

\noindent
{\bf I.5. The $K$-orbits of $S$.}
Under the action of $K$, the sphere $S=\{x\in V; Q(x,x)=1\}$
can be decomposed as a union of open and closed $K$-orbits.
The following Proposition is the special case of a more general
statement (see [MS, Prop. 3.9]).

\noindent
{\bf Proposition 1.7.} {\em Set $d_r={\rm diag}(\pii^r,1,1,1,\pii^{-r})$.
Each $K$-orbit of $S$ is of the form $Kd_rv_1$, $r\geq 0$.
The element ${}^t(x_1,x_2,x_3,x_4,x_5)\in S$ belongs to the $K$-orbit of
$d_rv_1$ if and only if
$$
\|x\|=\mbox{\rm max}\{|x_1|,|x_2|,|x_3|,|x_4|,|x_5|\}\ \mbox{is equal to}\ q^r.
$$ }
\noindent
{\em Proof.}
Consider the set $Kd_rv_1$. If ${\bf k}_i$ is the $i$th column of $k\in K$,
then $kd_rv_1$ is equal to
$\frac{1}{2}{\bf k}_1\pii^r+{\bf k}_5\pii^{-r}$.
When $k$ ranges over all elements of $K$, this sum ranges over
all elements ${}^t(x_1,x_2,x_3,x_4,x_5)\in S$ such that
$\mbox{max}\{|x_1|,|x_2|,|x_3|,|x_4|,|x_5|\}=q^r$, since the max absolute
value of the entries of $k_i$ is 1.
\hspace*{\fill}$\Box$

Let $\phi_r$ be the characteristic function of the $K$-orbit of
$d_rv_1$ in $S$.
Normalize the additive measure $dx$ on $F$ and the multiplicative measure
$d^\times x$ on $F^\times$ ($d^\times x=(1-1/q)^{-1}|x|^{-1}dx$), so that
$$
\int_Rdx=1,\qquad \mbox{and}\qquad \int_{|x|=1}d^\times x=
(1-1/q)^{-1}\int_{|x|=1}\frac{dx}{|x|}=1.
$$
Normalize the measure on $K$ so that its volume is $1$.
We need the following result.

\noindent
{\bf Proposition 1.8.} {\em The volume $\Lambda_r$ of the $K$-orbit of
$d_rv_1$ in $S$ is given by
\EQ\label{L2.1}
\Lambda_r=q^{3r}(1\mp q^{-2})\Lambda_0,\ \mbox{ if }\ r\geq 1,
\EN
where the ``$-$'' sign is in the split case and ``$+$'' in the non-split case. }

\noindent
{\em Proof.} Suppose $r\geq 1$.
Let $f_1$ be the characteristic function of $Kb_1K$ in $G$.
Since the measure $ds$ on $S$ is invariant under
the action of $G$, we have
\EQ\label{L2.2}
\int_G\int_S\phi_r(gs)dsf_1(g)dg=\int_Gf_1(g)dg\int_S\phi_r(s)ds.
\EN
Using Proposition 1.6, and that $\# N_i=q^{4-i}$,
the right hand side of (\ref{L2.2}) is equal to
\EQ\label{L2.3}
(q^3+q^2+q+1)\Lambda_r.
\EN
In the left hand side of (\ref{L2.2}), we change the order of integration.
It is
$$
\int_SI_r(s)ds,\ \mbox{ where }\   I_r(s)=\int_G\phi_r(gs)f_1(g)dg.
$$
Note that $I_r(ks)=I_r(s)$ for any $k\in K$. Thus it is constant on the
$K$ orbits in $S$. In particular, if $s$ is in the $K$-orbit of $d_iv_1$,
we have $I_r(s)=I_r(d_iv_1)$. We obtain
\EQ\label{L2.4}
\int_S I_r(s)ds=\sum_{i=0}^\infty \Lambda_iI_r(d_iv_1).
\EN
Using Proposition 1.6, the value of $I_r(d_iv_1)$ is
$$
\int_G\phi_r(gd_iv_1)f_1(g)dg  = \sum_{n\in N_1}\int_K\phi_r(kb_1nd_iv_1)dk
+\sum_{n\in N_2}\int_K\phi_r(kb_2nd_iv_1)dk
$$
$$
+ \sum_{n\in N_3}\int_K\phi_r(kb_2^{-1}nd_iv_1)dk
+\int_K\phi_r(kb_1^{-1}d_iv_1)dk.
$$
Since $\phi_r$ is left $K$-invariant, the expression above is equal to
\EQ\label{L2.4A}
\sum_{n\in N_1}\phi_r(b_1nd_iv_1)+\sum_{n\in N_2}\phi_r(b_2nd_iv_1)
+ \sum_{n\in N_3}\phi_r(b_2^{-1}nd_iv_1)+\phi_r(b_1^{-1}d_iv_1).
\EN
We consider the contribution to (\ref{L2.4}) from each sum of (\ref{L2.4A}).
In each case we distinguish between $r\geq 1$ and $r=0$. In some
cases we consider the case $r=1$ separately.
\np
{\em Case 1.} Consider the contribution to (\ref{L2.4}) from the first
sum in (\ref{L2.4A}). In this case $n\in N_1$ and the element $b_1nd_iv_1$,
in the split case, is equal to
$$
\left(\ColF{\pii^{-1}}{0}{0}{0}{0}\ColF{0}{1}{0}{0}{0}
\ColF{0}{0}{1}{0}{0}\ColF{0}{0}{0}{1}{0}\ColF{0}{0}{0}{0}{\pii}\right)
\left(\ColF{1}{0}{0}{0}{0}\ColF{-x_3}{1}{0}{0}{0}\ColF{-x_2}
{0}{1}{0}{0}\ColF{-x_1}{0}{0}{1}{0}\ColF{z}{x_1}{x_2}{x_3}{1}\right)
\left(\ColF{\frac{1}{2}\pii^i}{0}{0}{0}{\pii^{-i}}\right).
$$
In the non-split case the last vector column is
$$
{}^t(2\th\pii^i,0,0,1,\pii^{-i}).
$$
These elements are equal to (split/non-split cases respectively)
\EQ\label{L2.6}
{}^t(\pii^{i-1}/2+z\pii^{-(i+1)},x_1\pii^{-i},
x_2\pii^{-i},x_3\pii^{-i},\pii^{1-i}),
\EN
\EQ\label{V.6}
{}^t(2\th\pii^{i-1}-x_1\pii^{-1}+z\pii^{-(i+1)},x_1\pii^{-i},
x_2\pii^{-i},1+x_3\pii^{-i},\pii^{1-i}).
\EN
\indent
Depending on $n$, these elements belong to different $K$-orbits of $S$.
\np Let $r\geq 2$. Put $n(x_1,x_2,x_3)=n(x_1,x_2,x_3,0)$.
We decompose the group $N_1$ into a disjoint union as follows.
\np{\em (i)} Let $n=n(0,0,0)$ be the identity matrix.
Then $b_1nd_iv_1=b_1d_iv_1=d_{i-1}v_1$.
Hence $\phi_r(b_1nd_iv_1)=0$ unless $i-1=r$ (or $i=r+1$).
The contribution to (\ref{L2.4}) of this case is $\Lambda_{r+1}$.
\np{\em (ii)} Let $n=n(x_1,x_2,x_3)$ be a non-identity matrix with $z=0$.
We claim that there are $(q^2-1)$ such matrices. Indeed $z=0$ implies
$2x_1x_3+x_2^2=0$. The latter equation has $q(q-1)$ solutions with $x_3\neq 0$
($x_2$ arbitrary) and $q-1$ solutions with $x_2=x_3=0$, $x_1\neq 0$.
Since $z=0$, but $n$ is non-identity,
the element $b_1nd_iv_1$ belongs to the $K$-orbit of $d_iv_1$.
Hence $\phi_r(b_1nd_iv_1)=0$ unless $i=r$.
The contribution to (\ref{L2.4}) of this case is $(q^2-1)\Lambda_r$.
\np{\em (iii)} The remaining matrices $n=n(x_1,x_2,x_3)$ have $z\neq 0$.
Since the order of $N_1$ is $q^3$, there are $q^3-q^2$ such matrices.
The element $b_1nd_iv_1$ is in the $K$-orbit of $d_{i+1}v_1$ (since
$z\neq 0$).
The function $\phi_r(b_1nd_iv_1)$ is $0$ unless $i+1=r$ (or $i=r-1$).
Thus, the contribution to (\ref{L2.4}) of this case is
$(q^3-q^2)\Lambda_{r-1}$.
The contribution from the cases {\em (i)}, {\em (ii)} and {\em (iii)} is
(when $r\geq 2$)
\EQ\label{L2.5}
(q^3-q^2)\Lambda_{r-1}+(q^2-1)\Lambda_r+\Lambda_{r+1},
\EN
\indent
Let $r=1$. The only contributions to (\ref{L2.4}) occur when
$i=0,1,2$. If $i=2$, the element $b_1nd_2v_1$ (see (\ref{L2.6}))
is in the $K$-orbit of $d_1v_1$
precisely when $n$ is the identity matrix. The contribution to (\ref{L2.4}) is
$\Lambda_2$.
If $i=1$, the element $b_1nd_1v_1$ is in the $K$-orbit of $d_1v_1$ when
$z=0$ and $n$ is a non-identity matrix. As we have seen above,
the equation $z=0$ has $q^2-1$ solutions. The contribution is
$(q^2-1)\Lambda_1$. If $i=0$ the element $b_1nv_1$ is
\EQ\label{V.8}
{}^t((1/2+z)\pii^{-1},x_1,x_2,x_3,\pii),\
{}^t(\pii^{-1}(2\th-x_1+z),x_1,x_2,x_3,\th\pii)
\EN
(split/non-split cases respectively) is in the $K$-orbit of $d_1v_1$ when
$z+\frac{1}{2}\neq 0$ or $2\th-x_1+z\neq 0$.
The equation $z+\frac{1}{2}=0$
is equivalent to $2x_1x_3+x_2^2=1$, which has $q^2+q$ solutions.
Indeed, there are $q(q-1)$ solutions with $x_3\neq 0$, $x_2$ arbitrary,
and $2q$ solutions with $x_3=0$, $x_2=\pm 1$, $x_1$ arbitrary.
\np
The equation $2\th-x_1+z=0$ (where $z=-x_1x_3-x_2^2/2$)
is equivalent to $2x_1(1+x_3)+x_2^2=4\th$, which has $q^2-q$ solutions.
Indeed, there are $q(q-1)$ solutions with $x_1\neq 0$, $x_2$ arbitrary,
and no solutions with $x_1=0$, since $4\th$ is a non-square.
\np
Hence, this case contributes $(q^3-q^2-q)\Lambda_0$.
So, when $r=1$, we have
\EQ\label{L2.5a}
(q^3-q^2\mp q)\Lambda_0+(q^2-1)\Lambda_1+\Lambda_2,
\EN
where the ``$-$'' is in the split case and ``$+$'' in the non-split case.
\np
Let $r=0$. We distinguish between two cases:
\np{\em (i)} Let $n$ be the identity matrix, i.e.
$x_1=x_2=x_3=0$, and $z=0$. The element $b_1d_iv_1$ (see (\ref{L2.6}))
is in the $K$-orbit of $d_{i-1}v_1$. We have $\phi_0(d_{i-1}v_1)=0$
unless $i=1$. This contributes (in both split and non-split cases)
$\Lambda_1$ to (\ref{L2.4}).
\np{\em (ii)} Let $n$ be a non-identity matrix. Since
$x_1$, $x_2$, $x_3$ are not all zeroes, (\ref{L2.6}) is the $K$-orbit
of $d_{-i}v_1$. Hence $\phi_0(d_{-i}v_1)=0$ unless $i=0$.
Thus only $i=0$ contributes. This contribution occurs when
$b_1nv_1\in Kv_1$. This happens (see (\ref{V.8}))
when $1/2+z=0$, in the split case, or $2\th-x_1+z=0$ in the non-split case.
The first equation has $q^2+q$ solutions, the second one $q^2-q$.
Their contribution to (\ref{L2.4}) is $q(q\pm 1)\Lambda_0$.
\np
Thus Case 1, with $r=0$, contributes to (\ref{L2.4}) the quantity
$(q^2\pm q)\Lambda_0+\Lambda_1$.
\np
{\em Case 2.} Consider the contribution to (\ref{L2.4}) from the second
sum of (\ref{L2.4A}). For $n\in N_2$, the element $b_2nd_iv_1$, in the split
case is equal to
$$
\left(\ColF{1}{0}{0}{0}{0}\ColF{0}{\pii^{-1}}{0}{0}{0}
\ColF{0}{0}{1}{0}{0}\ColF{0}{0}{0}{\pii}{0}\ColF{0}{0}{0}{0}{1}\right)
\left(\ColF{1}{0}{0}{0}{0}\ColF{0}{1}{0}{0}{0}\ColF{0}
{-x_4}{1}{0}{0}\ColF{-x_1}{-\frac{1}{2}x_4^2}{x_4}{1}{0}
\ColF{0}{x_1}{0}{0}{1}\right)
\left(\ColF{\frac{1}{2}\pii^i}{0}{0}{0}{\pii^{-i}}\right)=
\left(\ColF{\frac{1}{2}\pii^i}{x_1\pii^{-(i+1)}}{0}{0}{\pii^{-i}}\right),
$$
and in the non-split case, replacing $(\pii^i/2,0,0,0,\pii^{-i})$
with $(2\th\pii^i,0,0,1,\pii^{-i})$, we obtain
$$
\left(2\th\pii^i-x_1,\pii^{-1}(x_1\pii^{-i}-\frac{1}{2}x_4^2),
x_4,\pii,\pii^{-i}\right).
$$
\indent
For $r\geq 1$ in the split case and $r\geq 2$ in the non-split case, we
have
\np{\em (i)} If $x_1=0$, the element $b_2nd_iv_1$
belongs to the $K$-orbit of $d_iv_1$. Since $x_4$ can be arbitrary,
the contribution to (\ref{L2.4}) is $q\Lambda_r$.
\np{\em (ii)} If $x_1\neq 0$ (there are $q^2-q$ such matrices),
the element $b_2nd_iv_1$ belongs to
the $K$-orbit of $d_{i+1}v_1$. We have $\phi_r(d_{i+1})=0$ unless
$i+1=r$ (or $i=r-1$). Their contribution is $(q^2-q)\Lambda_{r-1}$.
\np
Consider the non-split case with
$r=1$. The contribution occurs when $i=0,1$. If $i=1$, we have
that $x_1=0$ and $x_4$ can be arbitrary. If $i=0$, we have a
contribution when $x_1-\frac{1}{2}x_4^2\neq 0$, which has $q^2-q$
solutions. The resulting contribution is the same as in case $r\geq 2$.
\np
Thus, the contribution from this case (for $r\geq 1$) to (\ref{L2.4})
is (same in both split/non-split cases)
\EQ\label{L2.7}
(q^2-q)\Lambda_{r-1}+q\Lambda_r.
\EN
\indent
Now let $r=0$. First, consider the split case.
If $x_1\neq 0$, the element $b_2nd_iv_1$ lies
in the $K$-orbit of $d_{i+1}v_1$. Hence there are no positive $i$ for which
$b_2nd_iv_1$ belongs to the $K$-orbit of $v_1=d_0v_1$. If $x_1=0$ then
$b_2nd_iv_1$ is in the $K$-orbit of $d_iv_1$.
We have $\phi_0(d_iv_1)=0$ unless $i=0$. Since
$x_1$ is arbitrary, this case contributes $q\Lambda_0$ to (\ref{L2.4}).
In the non-split case, the only contribution occurs when $i=0$ and
$x_1-\frac{1}{2}x_4^2=0$. This case contributes $q\Lambda_0$ to (\ref{L2.4}).
\newline\indent
{\em Case 3.} Consider the contribution to (\ref{L2.4}) from the third
sum of (\ref{L2.4A}). For $n\in N_3$, the element $b_2^{-1}nd_iv_1$ is
$$
\left(\ColF{1}{0}{0}{0}{0}\ColF{0}{\pii}{0}{0}{0}
\ColF{0}{0}{1}{0}{0}\ColF{0}{0}{0}{\pii^{-1}}{0}\ColF{0}{0}{0}{0}{1}\right)
\left(\ColF{1}{0}{0}{0}{0}\ColF{-x_3}{1}{0}{0}{0}\ColF{0}
{0}{1}{0}{0}\ColF{0}{0}{0}{1}{0}\ColF{0}{0}{0}{x_3}{1}\right)
\left(\ColF{\frac{1}{2}\pii^i}{0}{0}{0}{\pii^{-i}}\right)=
\left(\ColF{\frac{1}{2}\pii^i}{0}{0}{x_3\pii^{-(i+1)}}{\pii^{-i}}\right)
$$
in the split case, and in the non-split case it is
$$
\left(2\th\pii^i,0,0,\pii^{-1}(1+x_3\pii^{-i}),\pii^{-i}\right).
$$
Consider $r\geq 1$. If $x_3=0$ the element $b_2nd_iv_1$
lies in the $K$-orbit of $d_iv_1$; otherwise it is in the $K$-orbit of
$d_{i+1}v_1$. Thus, the only contribution to (\ref{L2.4}) occurs when $i=r$
or $i=r-1$. Since $N_3$ has $q$ elements, this contribution
(for both split and non-split cases) is
\EQ\label{L2.8}
(q-1)\Lambda_{r-1}+\Lambda_r.
\EN
If $r=0$, the only contribution, $\Lambda_0$,
to (\ref{L2.4}), occurs when $i=0$ and $x_3=0$ or $x_3+1=0$ in the split
and non-split cases respectively.
\newline\indent
{\em Case 4.} The element $b_1^{-1}d_iv_1$ is in the $K$-orbit of
$d_{i+1}v_1$. When $r\geq 1$, it contributes the term $\Lambda_{r-1}$
to (\ref{L2.4}). There is no contribution when $r=0$.
\newline\indent
In all expressions below, the upper sign corresponds to the split case
and the lower to a non-split case.
\np
Summing up, when $r\geq 2$, the sum (\ref{L2.4}) is equal to the sum of
(\ref{L2.5}), (\ref{L2.7}), (\ref{L2.8}) and $\Lambda_{r-1}$:
$$
q^3\Lambda_{r-1}+(q^2+q)\Lambda_r+\Lambda_{r+1}.
$$
In all expressions below, the upper sign corresponds to the split case
and the lower to a non-split case.
When $r=0$, the sum (\ref{L2.4}) is equal to
$$
(q^2\pm q)\Lambda_0+\Lambda_1+q\Lambda_0+\Lambda_0
=(q^2+2q+1)\Lambda_0+\Lambda_1.
$$
When $r=1$, using (\ref{L2.5a}) instead of (\ref{L2.5}), this sum is
$$
(q^3\mp q)\Lambda_0+(q^2+q)\Lambda_1+\Lambda_2.
$$
\indent
Hence, we obtained the left hand sides of the equation (\ref{L2.2}),
for $r\geq 2$, $r=0$ and $r=1$.
Using (\ref{L2.3}), when $r=0,1$, we have
$$
\begin{array}{rcl}
(q^3\mp q)\Lambda_0+(q^2+q)\Lambda_1+\Lambda_2 & = & (q^3+q^2+q+1)\Lambda_1,
\\ \\ (q^2+q\pm q+1)\Lambda_0+\Lambda_1 & = & (q^3+q^2+q+1)\Lambda_0.
\end{array}
$$
The equations imply that
$$
\Lambda_1=(q^3\mp q)\Lambda_0=q^3(1\mp q^{-2})\Lambda_0,\ \mbox{and}\
\Lambda_2=q^6(1\mp q^{-2})\Lambda_0.
$$
For $r\geq 2$, we have
$$
q^3\Lambda_{r-1}+(q^2+q)\Lambda_r+\Lambda_{r+1} = (q^3+q^2+q+1)\Lambda_r,
$$
namely
$\Lambda_{r+1}-\Lambda_r=q^3(\Lambda_r-\Lambda_{r-1})$, or
$$
\Lambda_{r+1}=(1+q^3)\Lambda_r-q^3\Lambda_{r-1}.
$$
The proposition follows by induction.
\hspace*{\fill}$\Box$

\begin{center}
{\bf II. Correspondence of spherical functions.}
\end{center}

\noindent
{\bf II.1. Satake transform on the Hecke algebra of $G$.}
Let $\pi=I_G(\zeta,\zeta')$ be the unramified
representation of $G$ normalizedly induced from the character
$|\al_1|^\zeta\!|\al_2|^{\zeta'}$ of $B$.
Define the Satake transform $f^\vee$ of a spherical function $f$ on $G$
(i.e. $f\in C_c^\infty(G)$ and it is $K$-biinvariant), by $f^\vee(\pi)=
\mbox{tr}\ \pi(f)$. Then
$$
\pi(f)\phi_0=f^\vee(\pi)\phi_0,
$$
where $\phi_0$ is the unique $-$ up to a scalar multiple $-$
$K$-fixed vector in the space of $I_G(\zeta,\zeta')$.
We fix $\phi_0$ by $\phi_0(1)=1$.
\np
Set $\pi_\zeta=I_G(\zeta,1/2+\zeta)$.
We shall show in Proposition 2.2 below that there exists a unique
$K$-invariant function $T_\zeta$ on the sphere $S$, such that
$T_\zeta(1)\neq 0$ and
\EQ\label{40.1}
\int_Gf(g)T_\zeta(gs)dg=f^\vee(\pi_\zeta)T_\zeta(s).
\EN
Applying (\ref{40.1}) with $s=1$, the Satake transform of any spherical
function $f$ on $G$ is given on the $\pi_\zeta$ by
\EQ\label{40.1a}
f^\vee(\pi_\zeta)=\int_S\phi_f(s)T_\zeta(s)T_\zeta(1)^{-1}ds.
\EN
\indent
Since $T_\zeta$ is $K$-invariant, it will be defined by its values on
the $K$-orbits of $S$, which are of the form $Kd_rv_1$,
where $d_r=\dia(\pii^r,1,1,1,\pii^{-r})$ and $r\geq 0$.
Set $T_{r,\zeta}=T_\zeta(d_rv_1)$, $r\geq 0$.
These numbers are computed in Proposition 2.2.
\np
Put $\Phi_r=\sum_{i=0}^r\phi_i$, where $\phi_i$ is
the characteristic function of the $K$-orbit of $d_iv_1$.
Then, the function $\Phi_r$ is the
characteristic function of a subset of $S$ defined by
$$
\{ {}^t(x_1,x_2,x_3,x_4,x_5)\in S;\ |x_i|\leq q^r\!,\, 1\leq i\leq 5\}.
$$
The main goal of this subsection is to compute the integral
$\int_S \Phi_r(s)T_\zeta(s)ds$, $r\geq 1$.
\np
To compute the numbers $T_{r,\zeta}$, we need the following
result. Recall that $f_1$ is the characteristic
function of the double coset $Kb_1K$, where $b_1$ is
$\dia(\pii^{-1},1,1,1,\pii)$.

\noindent
{\bf Proposition 2.1.} {\em The Satake transform $f_1^\vee(\pi)$ at
$\pi=I_G(\zeta,\zeta')$ of $f_1$ is
$$
f_1^\vee(\pi)=q^{3/2}(q^{\zeta}+q^{\zeta'-\zeta}+q^{\zeta-\zeta'}+
q^{-\zeta}).
$$ }
\noindent
{\em Proof.} We follow [FM], Section F. Any
$\phi\in I_G(\zeta,\zeta')$ satisfies
\EQ\label{40.2}
\phi(nak)=\delta_B(a)^{1/2}|\al_1(a)|^\zeta|\al_2(a)|^{\zeta'}\phi(k),
\EN
where $a=\dia(\al,\beta,1,\beta^{-1},\al^{-1})$, and
$\al_1$ and $\al_2$ are the simple roots of $G=SO(5)$, defined by
$\al_1(a)=\al/\beta$ and $\al_2(a)=\beta$. We have
$$
(\pi(f)\phi)(h)=\int_Gf(g)\phi(hg)dg.
$$
Using the measure decomposition $dg=\delta_B(a)^{-1}dndadk$
and (\ref{40.2}), this is equal to
$$
\int_{N_0}\int_{A_0}\int_K f(h^{-1}nak)\delta_B(a)^{-1/2}|\al_1(a)|^\zeta
|\al_2(a)|^{\zeta'}\phi(k)dndadk.
$$
Put
$$
F_{f}(a)=\delta_B(a)^{-1/2}\int_K\int_{N_0} f(k^{-1}ank)dndk.
$$
Then
$$
\mbox{tr}\,I(\zeta,\zeta';f)=\int_{A_0}F_{f}(a)
|\al\beta^{-1}|^\zeta|\beta|^{\zeta'}da.
$$
When $f$ is $K$-biinvariant,
$$
F_f(a)=\delta_B(a)^{-1/2}\int_{N_0} f(an)dn.
$$
\indent
According to Proposition 1.6, the double coset $Kb_1K$ is the disjoint union
$$
Kb_1N_1\cup Kb_2N_2\cup Kb_2^{-1}N_3\cup Kb_1^{-1}.
$$
It follows that the integral $F_{f_1}(a)$ vanishes unless $a$ is
in the $K$-double cosets of
$b_1$, $b_2$, $b_2^{-1}$ or $b_1^{-1}$. Further,
$\delta_B(b_1)=q^3$, $\delta_B(b_2)=q$,
$\delta_B(b_2^{-1})=q^{-1}$ and  $\delta_B(b_1^{-1})=q^{-3}$. Hence
$$
\begin{array}{l}
F_{f_1}(b_1) = q^{-3/2}q^3 = q^{3/2},\
F_{f_1}(b_2)  = q^{-1/2}q^2 = q^{3/2}, \\\\
F_{f_1}(b_2^{-1}) = q^{1/2}q = q^{3/2},\
F_{f_1}(b_1^{-1}) = q^{3/2}1 = q^{3/2}.
\end{array}
$$
Evaluating the characters $\al_1$, $\al_2$ at $b_1$, $b_2$, $b_2^{-1}$
and $b_1^{-1}$, we obtain
$$
\int_{A_0}F_{f_1}(a)|\al\beta^{-1}|^\zeta|\beta|^{\zeta'}da=
q^{3/2}(q^\zeta+q^{\zeta'-\zeta}+q^{\zeta-\zeta'}+q^{-\zeta}).
$$
Since $f_1^\vee(\pi)=\ \mbox{tr}\ I_G(\zeta,\zeta';f_1)$, the
proposition follows.
\hspace*{\fill}$\Box$

We use Proposition 2.1 to prove the following:

\noindent
{\bf Proposition 2.2.} {\em The equation (\ref{40.1}) has a unique
solution satisfying $T_{0,\zeta}=1$, and (for $r\geq 1$), we have
$$
T_{r,\zeta}=\sum_{\xi\in\{\zeta,-\zeta\}}
\frac{(q^{\frac{3}{2}+\xi}-1)(1\mp q^{-\frac{1}{2}-\xi})}
{(q^{\xi}-q^{-\xi})(q^{\frac{3}{2}}\mp q^{-\frac{1}{2}})}
q^{-r(\frac{3}{2}-\xi)},
$$
where the ``$+$'' sign occurs in the split case and ``$-$'' in the
non-split case. }

\noindent
{\em Proof.} By Proposition 2.1, the equation (\ref{40.1}),
with $s=d_rv_0$ and $f=f_1$, becomes
\EQ\label{40.3}
\int_G f_1(g)T_\zeta(gd_rv_0)dg=q^{3/2}(q^{\zeta}+q^{1/2}+
q^{-1/2}+q^{-\zeta})T_{r,\zeta}.
\EN
Since $f_1$ is the characteristic function of the double coset
$Kb_1K=Kb_1N_1\cup Kb_2N_2\cup Kb_2^{-1}N_3\cup Kb_1^{-1}$,
and the function $T_\zeta$ is invariant under $K$, whose volume is $1$,
the left hand side of (\ref{40.3}) equals
\EQ\label{40.3a}
\sum_{n\in N_1}T_\zeta(b_1nd_rv_1)+\sum_{n\in N_2}T_\zeta(b_2nd_rv_1)+
\sum_{n\in N_3}T_\zeta(b_2^{-1}nd_rv_1)+T_\zeta(b_1^{-1}d_rv_1).
\EN
We will compute each of the sums.
\np
To simplify the notations (in the proof of this proposition), we write
$T$ for $T_\zeta$ and $T_r$ for $T_{r,\zeta}$.
For more details see Proposition 1.8.
\newline\indent
{\em Case 1.} Consider the first term of (\ref{40.3a}).
Put $n(x_1,x_2,x_3)=n(x_1,x_2,x_3,0)$.
\newline
Let $r\geq 1$. In the split case, we have
$$
b_1n(x_1,x_2,x_3)d_rv_1={}^t(\pii^{r-1}/2+z\pii^{-(r+1)},x_1\pii^{-r},
x_2\pii^{-r},x_3\pii^{-r},\pii^{1-r}).
$$
In the non-split case, we have
$$
b_1n(x_1,x_2,x_3)d_rv_1={}^t(2\th\pii^{r-1}-x_1\pii^{-1}+z\pii^{-(r+1)},
x_1\pii^{-r},x_2\pii^{-r},1+x_3\pii^{-r},\pii^{1-r}).
$$
Consider the following cases:
\np{\em (i)} Let $n=n(0,0,0)$ be the identity matrix.
Then $b_1nd_rv_1=b_1d_rv_1=d_{r-1}v_1$ and
$T(b_1nd_rv_1)=T(d_{r-1}v_1)=T_{r-1}$
\np{\em (ii)} Let $n=n(x_1,x_2,x_3)$ be a non-identity matrix with $z=0$.
In Proposition 1.8, we showed that there are $(q^2-1)$ such matrices.
The element $b_1nd_rv_1$ belongs to the $K$-orbit of $d_rv_1$.
Hence, $T(b_1nd_rv_1)=T(d_rv_1)=T_r$.
\np{\em (iii)} The remaining $q^3-q^2$ matrices $n=n(x_1,x_2,x_3)$ have $z\neq 0$.
The element $b_1nd_rv_1$ is in the $K$-orbit of $d_{r+1}v_1$.
Thus, $T(b_1nd_rv_1)=T(d_{r+1}v_1)=T_{r+1}$.
\np
We conclude that for $r\geq 1$, we have
\EQ\label{40.4}
\sum_{n\in N_1}T(b_1nd_rv_1)=(q^3-q^2)T_{r+1}+(q^2-1)T_r+T_{r-1}.
\EN
Let $r=0$ ($d_0=1$). The element $b_1n(x_1,x_2,x_3)v_1$
is equal to
$$
{}^t(\pii^{-1}(1/2+z),x_1,x_2,x_3,\pii),\ \mbox{or}\
{}^t(\pii^{-1}(2\th-x_1+z),x_1,x_2,1+x_3,\pii).
$$
in the split/non-split cases respectively.
These elements are in the $K$-orbit of $d_1v_1$ if $\frac{1}{2}+z\neq 0$
or $2\th-x_1+z\neq 0$. Otherwise, they are in the $K$ orbit of $v_1$.
\np
The equation $\frac{1}{2}+z=0$ has $q^2+q$ solutions, and
$2\th-x_1+z=0$ has $q^2-q$ solutions. Thus
\EQ\label{40.5}
\sum_{n\in N_1}T(b_1nv_1)=(q^3-q^2\mp q)T_1+(q^2\pm q)T_0.
\EN
\indent
{\em Case 2.} Consider the second term of (\ref{40.3a}).
The element $b_2n(x_1,0,0,x_4)v_1$
is equal to
$$
{}^t(\pii^r/2,x_1\pii^{-(r+1)},0,0,\pii^{-r}),\ \mbox{or}\
{}^t(2\th\pii^r-x_1,\pii^{-1}(x_1\pii^{-r}-x_4^2/2),x_4,\pii,\pii^{-r})
$$
in the split and non-split cases respectively.
\newline\noindent
Let $r\geq 0$ in the split cases or $r\geq 1$ in a non-split
case. We have:
\np{\em (i)} If $n=n(0,0,0,x_4)\in N_2$
the element $b_2nd_rv_1$ is in the $K$-orbit of $d_rv_1$.
Thus, $T(b_2nd_rv_1)=T(d_rv_1)=T_r$.
\np{\em (ii)} If $n=n(x_1,0,0,x_4)\in N_2$ has $x_1\neq 0$,
the element $b_2nd_rv_1$ is in the $K$-orbit of $d_{r+1}v_1$ and
$T(b_2nd_rv_1)=T(d_{r+1}v_1)=T_{r+1}$.
\np
Since there are $q$ matrices in {\em (i)} and $q^2-q$ in {\em (ii)},
(for $r\geq 0$), we have
\EQ\label{40.6}
\sum_{n\in N_2}T(b_2nd_rv_1)=(q^2-q)T_{r+1}+qT_r.
\EN
Let $r=0$. The element $b_2nv_1$ is in the $K$-orbit of $d_1v_1$
if $x_1-\frac{1}{2}x_4^2\neq 0$ ($q^2-q$ cases), otherwise it is in
the $K$-orbit of $v_1$. The contribution is the same as (\ref{40.6}).
\np
{\em Case 3.}  Consider the third term of (\ref{40.3a}).
\newline Let $r\geq 0$. The element $b_2^{-1}n(0,0,x_3,0)v_1$ is equal to
$$
{}^t(\pii^r/2,0,0,x_3\pii^{-(r+1)},\pii^{-r}),\ \mbox{or}\
{}^t(2\th\pii^r,0,0,\pii^{-1}(1+x_3\pii^{-r})\pii^{-r})
$$
in the split and non-split cases respectively.
This element belongs to the $K$-orbit of $d_rv_1$ if $n$ is the identity
matrix, and is in the $K$-orbit of $d_{r+1}v_1$ in the remaining $q-1$
cases. So, (for $r\geq 0$), we have
\EQ\label{40.8}
\sum_{n\in N_3}T(b_2^{-1}nd_rv_1)=(q-1)T_{r+1}+T_r.
\EN
\indent
{\em Case 4.} Since $b_1^{-1}d_r=d_{r+1}$ ($r\geq 0$), the last summand
of (\ref{40.3a}) is
\EQ\label{40.10}
T(b_1^{-1}d_rv_1)=T_{r+1}.
\EN
\indent
Adding (\ref{40.4}), (\ref{40.6}), (\ref{40.8}) and (\ref{40.10}) we
obtain that (when $r\geq 1$), the sum (\ref{40.3a}) is equal to
$$
q^3T_{r+1}+(q^2+q)T_r+T_{r-1}.
$$
When $r=0$, adding (\ref{40.5}), (\ref{40.6}), (\ref{40.8}) and
(\ref{40.10}) (used with $r=0$), the sum (\ref{40.3a}) is
$$
(q^3\mp q)T_1+(q^2+q\pm q+1)T_0,
$$
where the upper sign occurs in the split case and the lower in the non-split
case. Both expressions should be equal to the right hand side of (\ref{40.3}).
Thus, we obtained two difference equations
$$
q^3T_{r+1}+(q^2+q)T_r+T_{r-1}=q^{\frac{3}{2}}(q^\zeta+q^{\frac{1}{2}}+
q^{-\frac{1}{2}}+q^{-\zeta})T_r,
$$
and
$$
(q^3\mp q)T_1+(q^2+q\pm q+1)T_0=q^{\frac{3}{2}}(q^\zeta+q^{\frac{1}{2}}+
q^{-\frac{1}{2}}+q^{-\zeta})T_0.
$$
The first one can be simplified to
\EQ\label{40.11}
q^3T_{r+1}-q^{\frac{3}{2}}(q^\zeta+q^{-\zeta})T_r+T_{r-1}=0
\EN
In the second one we assume $T_0=1$. Then
\EQ\label{40.12}
T_1=\frac{q^{\frac{3}{2}}(q^\zeta+q^{-\zeta})\mp q-1}{q^3\mp q}.
\EN
The general solution of (\ref{40.11}) is given by
$$
T_r=c_1\lam_1^r+c_2\lam_2^r,
$$
where $\lam_1$ and $\lam_2$ are the two roots of the quadratic
equation
\EQ\label{40.13}
q^3\lam^2-q^{\frac{3}{2}}(q^\zeta+q^{-\zeta})\lam+1=0,
\EN
and $c_1$, $c_2$ are chosen to satisfy two initial conditions:
\EQ\label{40.14}
1=c_1+c_2,\mbox{ and } T_1=c_1\lam_1+c_2\lam_2.
\EN
The solutions of (\ref{40.13}) are $\lam_1=q^{-\frac{3}{2}+\zeta}$,
$\lam_2=q^{-\frac{3}{2}-\zeta}$ and that of (\ref{40.14}) are
$$
c_1=\frac{T_1-\lam_2}{\lam_1-\lam_2}=
\frac{(q^{\frac{3}{2}+\zeta}-1)(1\mp q^{-\frac{1}{2}-\zeta})}
{(q^\zeta-q^{-\zeta})(q^{\frac{3}{2}}\mp q^{-\frac{1}{2}})},
$$
$$
c_2=\frac{\lam_1-T_1}{\lam_1-\lam_2}=
-\frac{(q^{\frac{3}{2}-\zeta}-1)(1\mp q^{-\frac{1}{2}+\zeta})}
{(q^\zeta-q^{-\zeta})(q^{\frac{3}{2}}\mp q^{-\frac{1}{2}})}.
$$
The proposition follows.
\hspace*{\fill}$\Box$

The main result of this subsection is:

\noindent
{\bf Proposition 2.3.} {\em Set $X=q^\zeta$. Then
$$
\int_S\Phi_r(s)T_\zeta(s)ds=q^r\left[q^{\frac{1}{2}r}
\frac{X^{r+1}-X^{-(r+1)}}{X-X^{-1}}\mp q^{\frac{1}{2}(r-1)}
\frac{X^r-X^{-r}}{X-X^{-1}}\right],
$$
where the ``$-$'' sign is in the split case and the ``$+$'' in
the non-split case. }

\noindent
{\em Proof.} We have
\EQ\label{40.15}
\int_S\Phi_r(s)T(s)ds=\int_S\sum_{k=0}^r\phi_k(s)T(s)ds=\sum_{k=0}^rT_k
\Lambda_k.
\EN
Recall that $T_0=1$ and assume that $\Lambda_0=1$.
In the computations below, the upper sign occurs in the split case
and the lower in the non-split case.
Using Propositions 1.4 and 2.2,
we have that, the sum (\ref{40.15}) is equal to
$$
1+\frac{(q^{\frac{3}{2}}X-1)(1\mp q^{-\frac{1}{2}}X^{-1})}
{q^{\frac{3}{2}}(X-X^{-1})}\sum_{k=1}^rq^{\frac{3}{2}k}X^k-
\frac{(q^{\frac{3}{2}}X^{-1}-1)(1\mp q^{-\frac{1}{2}}X)}
{q^{\frac{3}{2}}(X-X^{-1})}\sum_{k=1}^rq^{\frac{3}{2}k}X^{-k}.
$$
Using the summation formula for geometric series, and then simplifying
the result we obtain the formula claimed in the proposition.
\hspace*{\fill}$\Box$

\noindent
{\bf II.2. The group H.} Recall that $H=PGL(2,F)$, and
$\pi_\zeta'=I_{H,\chi_0}(\zeta,-\zeta)$ denotes the representation
of $H$, induced from the character $\dia(a,1)n\mapsto |a|^\zeta\chi_0(a)$
of $B'$, where $\chi_0$ is
$1$ or the unramified character of $F^\times$, whose square is 1.
We denote the elements of $PGL(2,F)$ by their representatives in $GL(2,F)$.
The Bruhat decomposition of $H$ is $H=B'\cup N'wB'$, $B'=N'A'$, where
$A'=\{\dia(a,1);a\in F^\times\}$,
$$
N'=\left\{n(x)=\Two{1}{x}{0}{1};\  x\in F\right\},\
w=\Two{0}{1}{-1}{0}.
$$
The character $\psi'$ of $N'$ is defined by $\psi_{N'}(n(x))=\psi(x)$.
\newline\indent
Let $W_\zeta$ be the normalized unramified Whittaker function in the space of
representation $\pi_\zeta'$. It satisfies $W_\zeta(e)=1$,
$W_\zeta(ngk)=\psi_{N'}(n)W_\zeta(g)$
($n'\in N'$, $k\in K'=PGL(2,R)$, the standard maximal compact subgroup
of $H$), and for any $f'\in C_c^\infty(K'\bksl H/K')$, also
$$
f'^\vee(\pi_\zeta')W_\zeta(h)=\int_Hf'(x)W_\zeta(hx)dx.
$$
In particular, when $h=e$ we have
\EQ\label{40.1aa}
f'^\vee(\pi_\zeta')=\int_Hf'(x)W_\zeta(x)dx.
\EN
Since $dg=|a|^{-1}dndadk$, where $a=\dia(\al,1)$, $|a|=|\al|$ and
$da=d^\times\al$, we obtain
$$
f'^\vee(\pi_\zeta')=\int_H f'(nak)W_\zeta(nak)|a|^{-1}dndadk
$$
$$
=\int_{A'}\int_{N'}\!f'(na)\psi_{N'}(n)dn\,W_\zeta(a)|a|^{-1}da
=\int_{A'}\!\phi_{f'}'(a)W_\zeta(a)|a|^{-1}da,
$$
where
$$
\phi'_{f'}(g)=\int_{N'}f'(ng)\psi_{N'}(n)dn.
$$
The function $\phi'_{f'}$ lies in the space $\Ch$ of the right $K'$-invariant,
compactly supported modulo $N'$, functions $\phi'$ on $H$, which satisfy
$\phi'(ng)=\ovr\psi_{N'}(n)\phi'(g)$.
For any integer $r\geq 0$, define $\phi_r'$ by
$$
\phi_r'\left(\Two{\al}{0}{0}{1}\right)=\left\{\begin{array}{cc} 1, &\
\mbox{ if } |\al|=q^{-r}, \\ 0, &\ \mbox{ otherwise }.\end{array}\right.
$$

\noindent{\bf Proposition 2.4.}{\em
\np (1) The set $\{\phi'_r;r\geq 0\}$ is a basis of $\Ch$.
\np (2) For any $\phi'\in\Ch$, we have
$\phi'({\rm\dia}(\al,1))=0$ if $|\al|>1$.
\np (3) As in Proposition 2.3, put $X=q^{\zeta}$. Then
$$
\int_{A'}\phi_r'(a)W_\zeta(a)|a|^{-1}da
=(-1)^rq^{\frac{1}{2}r}\frac{X^{r+1}-X^{-(r+1)}}{X-X^{-1}}.
$$ }
\noindent
{\em Proof.} (1) This is clear.
\np (2) Indeed, choosing $n\in F^\times$ such that $|n|=|\al|$, we have
$$
\ovr\psi_{N'}(n)\phi'\left(\Two{\al}{0}{0}{1}\right)=
\phi'\left(\Two{1}{n}{0}{1}\Two{\al}{0}{0}{1}\right)=
\phi'\left(\Two{\al}{0}{0}{1}\Two{1}{\al/n}{0}{1}\right).
$$
This is $\phi'(\dia(\al,1))$, since $\phi'$ is right $K'$ invariant.
But $\psi_{N'}$ has conductor $R$, hence $\psi_{N'}(n)\neq 1$ for some $n$,
and our claim follows.
\np (3) This follows from the definition of $\phi_r'$ and Shintani's
explicit formula [Sh] for the Whittaker function (cf. [F], p. 305)
of $I_{H,\chi_0}(\zeta,-\zeta)$, which asserts that
$$
W_\zeta\left(\Two{\pii^r}{0}{0}{1}\right)=(-1)^rq^{-\frac{1}{2}r}
\frac{X^{r+1}-X^{-(r+1)}}{X-X^{-1}}.
$$
\hspace*{\fill}$\Box$

\noindent
{\bf II.3. The correspondence.}
Put $\pi_\zeta=I_G(\zeta,1/2+\zeta)$ and
$\pi_\zeta'=I_{H,\chi_0}{(\zeta,-\zeta)}$.
Following [FM], we say that $f\in C(K\bksl G/K)$ and $f'\in C(K'\bksl H/K')$
are {\it corresponding} if $f^\vee(\pi_\zeta)=f'^\vee(\pi_\zeta')$.
By (\ref{40.1a}) and (\ref{40.1aa}), an equivalent definition is given by
$$
\int_S \phi_f(s)T_\zeta(s)ds=\int_{A'}\phi_{f'}'(a)W_\zeta(a)|a|^{-1}da.
$$

\noindent
{\bf Definition.} Define a map
${\cal F}:C_c^\infty(K\bksl S)\rightarrow\Ch$
by ${\cal F}(\phi)=\phi'$ if
$$
\int_S \phi(s)T_\zeta(s)ds=\int_{A'}\phi'(a)W_\zeta(a)|a|^{-1}da,
$$
where $T_\zeta$ is the $K$-invariant function on $S$, defined in II.1
and Proposition 2.2, and $W_\zeta$ is the unramified normalized Whittaker
function defined in Section II.2.

\noindent
{\bf Proposition 2.5.} {\em The map ${\cal F}$ is well defined and induces a
linear bijection between the spaces $C_c^\infty(K\bksl S)$ and $\Ch$,
given by the correspondence ($r\geq 0$)
$$
{\cal F}(\Phi_r)=(-1)^rq^r(\phi_r'\pm \phi_{r-1}'),
$$
where the ``$+$'' sign occurs in the split case and the ``$-$'' sign in the
non-split case. }
\np
{\em Proof.} This follows from Propositions 2.3 and 2.4.
\hspace*{\fill}$\Box$

\noindent{\bf Corollary.} If $f$ and $f'$ are corresponding spherical functions,
then ${\cal F}(\phi_f)=\phi'_{f'}$.

\begin{center}
{\bf III. The Fourier coefficients of orbital integrals.}
\end{center}

\noindent
{\bf III.1. The Fourier coefficients of orbital integrals on $H$.}
As usual, $F$ is a $p$-adic field, and $\chi$ is a complex valued
character of $F^\times$ with conductor $m$, $m\geq 0$. Thus if $m\geq 1$,
this character is trivial on $1+\pii^mR$ and is non-trivial on $1+\pii^{m-1}R$.
If $m=0$ then $\chi$ is trivial on $R^\times$ and is non-trivial on
$1+\pii^{-1}R$, and we say that $\chi$ is unramified.
Recall that $\{\phi_r';r\geq 0\}$ is the basis of $\Ch$, and
for any $\phi'\in\Ch$, we defined
\EQ\label{6.1}
\Psi'(\al,\phi')=\int_{F^\times}
\phi'\left(\Two{0}{1}{-1}{0}\Two{1}{\al}{0}{1}\Two{a}{0}{0}{1}\right)
\chi_0(a)d^\times a.
\EN

\noindent
{\bf Definition.} For any $r\geq 0$, set
\EQ\label{6.2}
\widehat{\Psi}'_r(\chi)=\int_{F^\times}(\th,\al)\psi(\al)\chi(\al)|\al|
\Psi'(\al^{-1},\phi_r')d^\times\al.
\EN

\noindent
{\bf Proposition 3.1.} {\em Let $\chi$ be a multiplicative character of
$F^\times$. If $\chi$ is unramified, define $\zeta$ and $X$ by
$X=|\pii|^{-\zeta}=\chi(\pii)^{-1}$. Then
$$
\int_{F^\times}(\th,\al)\psi(\al)\chi(\al)|\al|\Psi'(\al^{-1}\!,(-1)^r
q^r(\phi_r'\pm \phi_{r-1}'))
d^\times\al
$$
is equal to $0$ if $\chi$ is ramified, and is equal to
$$
X^{-r}\pm qX^{1-r}
$$
if $\chi$ is unramified.}

This proposition follows from the following Proposition :

\noindent
{\bf Proposition 3.2.} {\em The integral
$$
\widehat{\Psi}_r'(\chi)=\int_{F^\times}(\th,\al)\psi(\al)\chi(\al)|\al|
\Psi'(\al^{-1}\!,\phi_r')d^\times\al
$$
is equal to
\EQ\label{8.1}
(-1)^r(qX)^{-r}\frac{1\pm qX}{1\mp qX} + (\mp 1)^{r-1}
\frac{2q^2X(1\mp X)}{(q-1)(1\mp qX)},
\EN
in the split and the non-split cases respectively,
if $\chi$ is unramified, and to
$$
2(-1)^rq^m\tau(\psi,\chi),\  \ \mbox{where}\ \
\tau(\psi,\chi)=\int_{|x|=q^m}\psi(x)\chi(x)d^\times x
$$
if $\chi$ is ramified with conductor $m$.}

Let us show how Proposition 3.1 follows from Proposition 3.2:

\noindent
{\em Proof of Proposition 3.1.}
If $\chi$ is ramified, the result is obvious. If $\chi$ is unramified,
using the result of Proposition 3.2, we have that
$$
(-1)^rq^r(\widehat{\Psi}'_r(\chi)\pm \widehat{\Psi}'_{r-1}(\chi))
$$
is equal to (in the split/non-split cases)
$$
X^{-r}\frac{1\pm qX\mp (qX\pm (qX)^2)}{1\mp qX}
=\frac{1-(qX)^2}{1\mp qX}X^{-r}=X^{-r}\pm qX^{1-r}.
$$
\hspace*{\fill}$\Box$

To prove Proposition 3.2, we need the following self contained Lemmas:

\noindent
{\bf Lemma 3.3.} {\em Let $\psi$ be a character of $F$ with conductor $R$
and $\chi$ an unramified character of $F^\times$,
where $F$ is $p$-adic field. Set $|\pii|=q^{-1}$ and $\chi(\pii)^{-1}=X$.
Then for any $x\in F^\times$,
we have
\EQ\label{6.6}
\int_{|\al|=q^k}|\al|^3\psi(\al x)\chi(\al)d^\times\al=\left\{
\begin{array}{ll} 0, & \mbox{ if }\ q^k\geq q^2|x|^{-1}, \\
-(q-1)^{-1}(q^3X)^k, & \mbox{ if }\ q^k=q|x|^{-1}, \\
(q^3X)^k, & \mbox{ if }\ q^k\leq |x|^{-1}.
\end{array}\right.
\EN }

\noindent{\em Proof.} Recalling that $d^\times\al=(1-1/q)^{-1}|\al|^{-1}d\al$,
and that if $|\al|=q^k$ then $\chi(\al)=X^k$, we obtain
$$
\int_{|\al|=q^k}|\al|^3\psi(\al x)\chi(\al)d^\times\al=
\cir^{-1}q^{2k}X^k\int_{|\al|=q^k}\psi(\al x)d\al.
$$
Put $\beta=\al x$. The latter integral is
$$
\cir^{-1}q^{2k}X^k|x|^{-1}\int_{|\beta|=q^k|x|}\psi(\beta)d\beta.
$$
Recall that
$$
\int_{|\beta|=q^l}\psi(\beta)d\beta=\left\{\begin{array}{ll}
0, & \mbox{ if }\  l\geq 2, \\
-1, & \mbox{ if }\ l=1, \\
(1-1/q)q^l, & \mbox{ if }\ l\leq 0.
\end{array}\right.
$$
The lemma follows from this if we note that if $l=1$, i.e. $q^k|x|=q$, then
$$
q^{2k}X^k|x|^{-1}=\frac{1}{q}(q^3X)^k.
$$
The other cases are obvious.
\hspace*{\fill}$\Box$

\noindent
{\bf Lemma 3.4.} {\em The same notations as in Lemma 3.3, but
let $\chi$ be a ramified character with conductor $m$, $m\geq 1$. Then
$$
\int_{|x|=q^k}\psi(x)\chi(x)d^\times x
$$
is equal to $0$ unless $k=m$, in which case we denote it by
$\tau(\psi,\chi)$.}

\noindent{\em Proof.} We will consider two cases: $k>m$ and $k<m$.
\np
{\em Case 1.} Consider the case $k>m$. In this case, there exists an
element $y\in F^\times$, such that $|y|>1$, $1+yx^{-1}\in 1+\pii^m R$
and $\psi(y)\neq 1$. For such $y$, we have
$$
\chi(1+y/x)=1,\ |x+y|=|x|.
$$
The change of variables $x\mapsto x+y$, gives
$$
\int_{|x|=q^k}\psi(x)\chi(x)d^\times x
=\int_{|x|=q^k}\psi(x+y)\chi(x(1+y/x))d^\times x
=\psi(y)\int_{|x|=q^k}\psi(x)\chi(x)d^\times x.
$$
This is $0$ since $\psi(y)\neq 1$.
\np
{\em Case 2.} Consider the case $k<m$. In this case, take an
element $y\in 1+\pii^mR$ (thus $\psi(y)=1$, $|xy|=|x|$) such that
$\chi(y)\neq 1$. Set $y'=y-1$. Since $|xy'|\leq 1$, we have
$$
\psi(xy)=\psi(x+xy')=\psi(xy')\psi(x)=\psi(x).
$$
The change of variables $x\mapsto xy$, gives
$$
\int_{|x|=q^k}\psi(x)\chi(x)d^\times x
=\int_{|x|=q^k}\psi(xy)\chi(xy)d^\times x
=\chi(y)\int_{|x|=q^k}\psi(x)\chi(x)d^\times x.
$$
This is $0$, since $\chi(y)\neq 1$. The lemma follows.
\hspace*{\fill}$\Box$

\noindent
{\em Proof of Proposition 3.2.} The integral $\widehat{\Psi}_r'(\chi)$ is given
by the integral
\EQ\label{8.2}
\int_{F^\times}\int_{F^\times}
\phi_r'\left(\Two{0}{1}{-1}{0}\Two{1}{\al^{-1}}{0}{1}\Two{a}{0}{0}{1}\right)
(\th,\al)\chi(\al)\iota(a)\psi(\al)|\al|d^\times ad^\times\al,
\EN
where $\iota(a)=\chi_0(a)$ in the split case and is $1$ in the non-split
case. By matrix multiplication
\EQ
\Two{1}{\al^{-1}}{0}{1}\Two{a}{0}{0}{1}=\Two{a}{0}{0}{1}\Two{1}{\al^{-1}a^{-1}}
{0}{1}.
\EN
We will consider two cases.
\newline\indent
{\em Case 1.} Assume that $|\al^{-1} a^{-1}|\leq 1$ or $|\al|\geq |a|^{-1}$.
Using (\ref{8.2}) and the property of $\phi_r'$
$$
\phi_r'\left(\!\Two{0}{1}{-1}{0}\!\!\Two{1}{\al^{-1}}{0}{1}\!\!\Two{a}{0}{0}{1}
\!\right)=\phi_r'\left(\!\Two{0}{1}{-1}{0}\!\!\Two{a}{0}{0}{1}\!\right)
=\phi_r'\left(\!\Two{1}{0}{0}{a}\!\right).
$$
In $PGL(2)$ this is equal to $\phi_r'\left(\Two{a^{-1}}{0}{0}{1}\right)$,
which is $0$ unless $|a|^{-1}=q^{-r}$ or $|a|=q^r$.
Hence, the integral (\ref{8.2}) is equal to
\EQ\label{8.4}
\int_{|a|=q^r}\iota(a)d^\times a\int_{|\al|\geq q^{-r}}(\th,\al)
\chi(\al)\psi(\al)|\al|d^\times\al=\iota^r(\pii)\sum_{l\geq -r}
\int_{|\al|=q^l}(\th,\al)\chi(\al)\psi(\al)|\al|d^\times\al.
\EN
If $\chi$ is unramified, we apply Lemma 3.3 with $|\al|$ instead of $|\al|^3$
and $x=1$:
\EQ\label{8.4a}
\int_{|\al|=q^l}\chi(\al)\psi(\al)
|\al|d^\times\al=\left\{\begin{array}{ll}  0, & \mbox{ if }\ l\geq 2,\\
-X(1-1/q)^{-1}, & \mbox{ if } l=1, \\
(qX)^{l}, &  \mbox{ if } l\leq 0. \end{array}\right.
\EN
Futhermore, in the split case $(\th,\al)=1$, and in the non-split case
$(\th,\al)=(-1)^l$ if $|\al|=q^l$. So, the sum (\ref{8.4}) is equal to
\EQ\label{8.5}
(\mp 1)^r\sum_{l=-r}^{0}(\pm qX)^l+(\mp 1)^{r-1}\frac{qX}{q-1}=
\frac{(-1)^r(qX)^{-r}+(\mp 1)^{r-1}qX}{1\mp qX}+(\mp 1)^{r-1}\frac{qX}{q-1}.
\EN
If $\chi$ is ramified with conductor $m$, then according to Lemma 3.4,
the integral (\ref{8.4}) is equal
to
$$
(-1)^rq^m\int_{|\al|=q^m}\chi(\al)\psi(\al)d^\times\al=
(-1)^rq^m\tau(\psi,\chi).
$$
\indent
{\em Case 2.} In this case, $|\al|< |a|^{-1}$. Note that
$$
\Two{0}{1}{-1}{0}\Two{1}{\al^{-1}}{0}{1}\Two{a}{0}{0}{1}
=\Two{0}{1}{-a}{-\al^{-1}}.
$$
In $PGL(2)$, we have
$$
\Two{0}{1}{-a}{-\al^{-1}}\Two{-\al}{0}{0}{-\al}=\Two{0}{-\al}{a\al}{1}
=\Two{1}{-\al}{0}{1}\Two{a\al^2}{0}{0}{1}
\Two{1}{0}{a\al}{1}.
$$
Hence,
$$
\phi_r'\left(\Two{0}{1}{-1}{0}\Two{1}{\al^{-1}}{0}{1}\Two{a}{0}{0}{1}
\right)=\psi(\al)\phi_r'\left(\Two{a\al^2}{0}{0}{1}\right).
$$
The integral (\ref{8.2}) reduces to
$$
\int_{F^\times}\int_{F^\times}\psi(\al)\phi_r'\left(\Two{a\al^2}{0}{0}{1}
\right)\iota(a)(\th,\al)\chi(\al)\psi(\al)|\al|d^\times\al d^\times a.
$$
This integral does not vanish precisely when $|a\al^2|=q^{-r}$. Thus, it is
\EQ\label{8.6}
\int_{|\al|< |a|^{-1}}\int_{|a|=q^{-r}|\al|^{-2}}
\iota(a)(\th,\al)\chi(\al)\psi(\al)|\al|d^\times\al d^\times a.
\EN
Set $|\al|=q^l$ and $|a|=q^s$. The above integral is taken over the set
$l<-s$, $s+2l=-r$. Equivalently, this set is defined by $l>-r$, $s=-r-2l$.
Applying (\ref{8.4a}) to (\ref{8.6}), the integral is a finite sum
(split/non-split cases respectively)
$$
\iota^r(\pii)\sum_{l>-r}\int_{|\al|=q^l}(\th,\al)\chi(\al)
\psi(\al)|\al|d^\times\al=
(\mp 1)^r\sum_{l=1-r}^0(\pm qX)^l+(\mp 1)^{r-1}\frac{qX}{q-1}
$$
if $\chi$ is unramified, and to $(-1)^rq^m\tau(\psi,\chi)$ if $\chi$
is ramified.
Using the summation formula for geometric series, the unramified case
is
$$
(\mp 1)^r\frac{(\pm qX)^{-r}-1}{(\pm qX)^{-1}-1}
+(\mp 1)^{r-1}\frac{qX}{q-1}=
\frac{(\mp 1)^r(\pm qX)^{1-r}\mp qX}{1\mp qX}+(\mp 1)^{r-1}\frac{qX}{q-1}.
$$
Combining this expression with (\ref{8.5}), we obtain the final result
(for the unramified case)
$$
\frac{(-1)^r(qX)^{-r}+(\mp 1)^{r-1}qX}{1\mp qX}+\frac{(-1)^{r-1}(qX)^{1-r}+
(\mp 1)^{r-1}qX}{1\mp qX}-\frac{2qX}{q-1}.
$$
Once simplified, it completes the proof of the Proposition 3.2.
\hspace*{\fill}$\Box$

\noindent
{\bf III.2. The Fourier coefficients of orbital integrals on $G$.}
Recall that for any $\phi\in C_c^\infty(K\bksl S)$ and
$a_\al=\dia(\al,1,1,1,\al^{-1})$, we defined
$$
\Psi(\al,\phi)=\int_N\phi(na_\al\gamma_0v_0)\psi_N(n)dn.
$$
For any $\phi\in C_c(K\bksl S)$, the Fourier transform
$\widehat{\Psi}_\phi(\chi)$ of $\Psi(\al,\phi)$ is given by
$$
\widehat{\Psi}_\phi(\chi)=\int_{F^\times}\Psi(\al,\phi)\chi(\al)d^\times\al.
$$
Recall that $\Phi_r=\sum_{i=0}^r\phi_i$, where $\phi_i$ is the
characteristic function of the $K$-orbit of $d_iv_1$.
In this section we compute
$$
\widehat{\Psi}_r(\chi)=\widehat{\Psi}_{\Phi_r}(\chi),
$$
where $\chi$ is the same character as in Section II.

\noindent{\bf Proposition 3.5.} {\em In the split case, the Fourier
transform $\widehat{\Psi}_r(\chi)$ of $\Psi(\al,\Phi_r)$ is equal to
$$
qX^{1-r}+X^{-r}
$$
if the character $\chi$ is unramified, and to $0$ if $\chi$ is ramified.}

\noindent{\em Proof.}
Recalling the definitions of $a$, $n$, $\gamma_0$ and $v_0$,
we have
$$
\Psi(\al,\Phi_r)=\int\!\!\!\int\!\!\!\int_{F^3}\Phi_r(\al/2+z/\al,x_1/\al,
x_2/\al,x_3/\al,1/\al)\psi(x_2)dx_1dx_2dx_3.
$$
Recall that $\Phi_r(x_1,x_2,x_3,x_4,x_5)$ is $1$ if $|x_i|\leq q^r$
$(i=1,...,5)$ and is zero otherwise. Thus, the integral
$\widehat{\Psi}_r(\chi)$ is equal to
\EQ\label{7.2}
\int\!\!\!\int\!\!\!\int\!\!\!\int
\chi(\al)\psi(x_2)dx_1dx_2dx_3d^\times\al,
\EN
over the set defined by
$$
\left|\frac{\al}{2}+\frac{z}{\al}\right|\leq q^r,\
\left|\frac{x_1}{\al}\right|\leq q^r\!,\
\left|\frac{x_2}{\al}\right|\leq q^r\!,\
\left|\frac{x_3}{\al}\right|\leq q^r\!,\
\left|\frac{1}{\al}\right|\leq q^r\!,
$$
where $z=-x_1x_3-x_2^2/2$.
After the change of variables $x_i\mapsto \al x_i$, the integral
(\ref{7.2}) is equal to
\EQ\label{7.4}
\int\!\!\!\int\!\!\!\int\!\!\!\int
|\al|^3\chi(\al)\psi(\al x_2)dx_1dx_2dx_3d^\times\al,
\EN
over the set
$$
|\al(1-x_2^2-2x_1x_3)|\leq q^r,\ |x_1|\leq q^r,\
|x_2|\leq q^r,\ |x_3|\leq q^r,\ q^{-r}\leq |\al|.
$$
Fixing $x_1$, $x_2$ and $x_3$, we obtain that
$$
q^{-r}\leq |\al|\leq q^r|1-x_2^2-2x_1x_3|^{-1}.
$$
Changing the order of integration in (\ref{7.4}), it is equal to
\EQ\label{7.7}
\int_{|x_1|\leq q^r}\int_{|x_3|\leq q^r}\int_{|x_2|\leq q^r}
\int_{q^{-r}\leq |\al|\leq q^r|1-x_2^2-2x_1x_3|^{-1}}
|\al|^3\chi(\al)\psi(\al x_2)d^\times\al dx_2dx_3dx_1.
\EN
First, assume that $\chi$ is ramified with conductor $m$. By Lemma 3.4
$$
\int_{|\al|=q^k}\chi(\al)\psi(\al)d^\times\al
$$
vanishes unless $k=m$. Thus
$$
\int_{|\al|=q^k}|\al|^3\chi(\al)\psi(\al x_2)d^\times\al=
q^{3m}\tau(\psi,\chi)|x_2|^{-3}\chi^{-1}(x_2).
$$
Since, for any $l$, the integral of $\chi^{-1}(x_2)$ over $|x_2|=q^l$ is
equal to $0$, we conclude that $\widehat{\Psi}_r(\chi)=0$ when $\chi$ is
ramified.
\newline\indent
Now, let us consider the case of an unramified character $\chi$.
To apply Lemma 3.3, we need to split the domain of integration
over $\al$ into two domains: defined by condition $q^r/|1-2x_1x_3-x_2^2|$
is $\leq 1/|x_2|$  or $> 1/|x_2|$. We will consider the contributions
of the integral (\ref{7.7}) over each of these domains. There are two cases.
\newline\indent
{\em Case 1.} Let us consider the contribution to (\ref{7.7}) from the
first domain, namely $q^r/|1-2x_1x_3-x_2^2|\leq 1/|x_2|$. Equivalently, it is
\EQ\label{7.8}
\left|\frac{1-2x_1x_3}{x_2}-x_2\right|\geq q^r.
\EN
\indent
In this case there are two subdomains.
\np{\em Case 1a.} The first subdomain is $|x_2|=q^r$. Since $|x_i|\leq q^r$,
we have $|1-2x_1x_3|\leq q^{2r}$, which implies
$$
\left|\frac{1-2x_1x_3}{x_2}-x_2\right|=q^r.
$$
This is equivalent to $|1-2x_1x_3-x_2^2|=|x_2|q^r=q^{2r}$.
Note that since $r\geq 1$, we have $|1-2x_1x_3-x_2^2|=|2x_1x_3+x_2^2|$.
Hence, we conclude that in the
integral (\ref{7.7}), the integration over $\al$ is taken over
$|\al|=q^{-r}$. Applying Lemma 3.3,
$$
\int_{|\al|=q^{-r}} |\al|^3\chi(\al)\psi(\al x_2)d^\times\al=(q^3X)^{-r}.
$$
The integral (\ref{7.7}) is the product of $(q^3X)^{-r}$ and the
volume of a subset defined by
$$
\{|x_1|\leq q^r\!\!,|x_2|\leq q^r\!\!,|x_3|\leq q^r\!\!,|2x_1x_3+x_2^2|=q^{2r}\}.
$$
This subset is equal to
\EQ\label{7.10}
\{|x_1x_3|\leq q^{2r-1}\!\!\!,\,|x_2|=q^r\}\cup
\{|x_1|=q^r\!\!,\,|x_3|=q^r\!\!,\,|x_2|=q^r\!\!,\,|2x_1x_3+x_2^2|=q^{2r}\}.
\EN
The volume of the first subset is
$$
\left[\int_{|x_1|<q^r}\int_{|x_3|\leq q^r}+
\int_{|x_1|=q^r}\int_{|x_3|<q^r}\right]dx_3dx_1\int_{|x_2|=q^r}dx_2
=\biggl(1-\frac{1}{q}\biggr)\biggl(2-\frac{1}{q}\biggr)q^{3r-1}.
$$
The volume of the second subset of (\ref{7.10}) is the integral
$$
\int_{|x_1|=q^r}\int_{|x_2|=q^r}
\int_{|x_3|=q^r,\, |2x_1x_3+x_2^2|=q^{2r}}dx_3dx_2dx_1
=\cir^2\biggl(1-\frac{2}{q}\biggr)q^{3r}.
$$
Multiplying the volume of (\ref{7.10}) by $(q^3X)^{-r}$,
the contribution of (\ref{7.7}) from the subcase 1a is:
\EQ\label{7.11}
\biggl(1-\frac{1}{q}\biggr)\biggl(2-\frac{1}{q}\biggr)\frac{1}{q}X^{-r}
+\cir^2\biggl(1-\frac{2}{q}\biggr)X^{-r}.
\EN
\newline\noindent
{\em Case 1b.} The second subdomain is defined by $|x_2|<q^r$. In this case,
we have
$$
\left|\frac{1-2x_1x_3}{x_2}-x_2\right|=\left|\frac{1-2x_1x_3}{x_2}\right|.
$$
Hence, in the integral (\ref{7.7}), the integration over $\al$
is performed over
$$
q^{-r}\leq |\al|\leq \frac{q^r}{|1-2x_1x_3|},
$$
and $x_2$ satisfies $\{ |x_2|<q^{r}\!\!,\,|x_2|\leq |1-2x_1x_3|q^{-r}\}$.
We will consider two cases: $|1-2x_1x_3|=q^{2r}$ and $|1-2x_1x_3|<q^{2r}$.
\newline\indent
{\em (i)} Let $|1-2x_1x_3|=q^{2r}$. Since $|x_i|\leq q^r$, this implies
that $|x_1|=|x_3|=q^r$. The integration (in (\ref{7.7})) is taken only over
$\al$ with $|\al|=q^{-r}$, and over $x_2$ with $|x_2|<q^r$.
Thus the integral (\ref{7.7}) is
$$
\int_{|x_1|=q^r}\int_{|x_3|=q^r}\int_{|x_2|<q^r}\int_{|\al|=q^{-r}}
|\al|^3\chi(\al)\psi(\al x_2)d^\times\al dx_2dx_3dx_1.
$$
Once evaluated and simplified, it is
\EQ\label{7.12}
\cir^2q^{2r}q^r\frac{1}{q}(q^3X)^{-r}=\cir^2\frac{1}{q}X^{-r}.
\EN
\indent
{\em (ii)} Let $|1-2x_1x_3|<q^{2r}$.
Define $l$ by $|1-2x_1x_3|=q^l$. In (\ref{7.7}), $x_2$ is bounded
from above by $|1-2x_1x_3|q^{-r}=q^{l-r}$. Since $l\leq 2r-1$, this
implies that $x_2<q^r$. The integral (\ref{7.7}) becomes
$$
\int_{|x_1|\leq q^r}\int_{|x_3|\leq q^r}
\int_{|x_2|\leq |1-2x_1x_3|q^{-r}}
\int_{q^{-r}\leq |\al|\leq q^{r}|1-2x_1x_3|^{-1}}
|\al|^3\chi(\al)\psi(\al x_2)d^\times\al dx_2dx_3dx_1.
$$
Breaking it into the sum over $l$, it is
\EQ\label{7.13}
\sum_{l\leq 2r-1}\int\!\!\int_{|x_1|\leq q^r,|x_3|\leq q^r,\ |1-2x_1x_3|=q^l}
\int_{|x_2|\leq q^{l-r}}\int_{q^{-r}\leq |\al|\leq q^{r-l}}
|\al|^3\chi(\al)\psi(\al x_2)d^\times\al dx_2dx_1dx_3.
\EN
Applying Lemma 3.3, the integral over $\al$ becomes a geometric series
$\sum_{k=-r}^{r-l}q^{3k}X^k$.
Substituting this into the integral over $x_2$ in (\ref{7.13}), we obtain
$$
\int_{|x_2|\leq q^{l-r}}\frac{(q^3X)^{r-l+1}-(q^3X)^{-r}}{q^3X-1}
dx_2=\frac{q(q^2X)^{r+1}}{q^3X-1}(q^2X)^{-l}-\frac{(q^4X)^{-r}}{q^3X-1}q^l.
$$
Hence, the integral (\ref{7.13}) is equal to
$$
\sum_{l\leq 2r-1}\int\!\!\int_{|x_1|\leq q^r,\,|x_3|\leq q^r,\,
|1-2x_1x_3|=q^l}\left[\frac{q(q^2X)^{r+1}}{q^3X-1}(q^2X)^{-l}-
\frac{(q^4X)^{-r}}{q^3X-1}q^l\right]dx_1dx_3.
$$
Splitting off the term corresponding to $l=2r-1$, it is
\EQ\label{7.15a}
\frac{1}{q}(q^3X+1)(q^2X)^{-r}\int\!\!\int_{|x_1|\leq q^r\!,\,
|x_3|\leq q^r\!,\,|1-2x_1x_3|=q^{2r-1}}dx_1dx_3
\EN
\EQ\label{7.15b}
+\sum_{l\leq 2r-2}\left[\frac{q(q^2X)^{r+1}}{q^3X-1}(q^2X)^{-l}-
\frac{(q^4X)^{-r}}{q^3X-1}q^l\right]\int\!\!\int_{|x_1|\leq q^r\!,\,
|x_3|\leq q^r\!,\,|1-2x_1x_3|=q^l}dx_1dx_3.
\EN
This is a contribution of (\ref{7.7}) over the domain 1b(ii).
We will not evaluate the sum (\ref{7.15b}) any further for it will be cancelled
by a similar sum obtained below.
\newline\indent
{\em Case 2.} Consider the contribution of (\ref{7.7}) over the domain
$q^r/|1-2x_1x_3-x_2^2|>1/|x_2|$. Equivalently, it is
\EQ\label{7.16}
\left|\frac{1-2x_1x_3}{x_2}-x_2\right|<q^r.
\EN
Using Lemma 3.3, the integral over $\al$ in (\ref{7.7}) can be
split into two integrals: one over the subdomain
$q^{-r}\leq |\al|\leq |x_2|^{-1}$
and the other one over the subdomain $|\al|=q|x_2|^{-1}$. Thus the integral
(\ref{7.7}) is
\EQ\label{7.17}
\int\!\!\int\!\!\int\left[\int_{q^{-r}\leq |\al|\leq |x_2|^{-1}}\!\!
+\int_{|\al|=q|x_2|^{-1}}\right]|\al|^3\psi(\al x_2)\chi(\al)d^\times\al
dx_1dx_2dx_3,
\EN
where $x_1$, $x_2$ and $x_3$ range over the set
\EQ\label{7.18}
|x_1|\leq q^r\!\!,\ |x_3|\leq q^r\!\!,\
\left|\frac{1-2x_1x_3}{x_2}-x_2\right|<q^r.
\EN
Define $l$ by $|1-2x_1x_3|=q^l$. We distinguish between the cases
$l\leq 2r-1$ and $l=2r$.
\newline\indent
{\em Case 2a}. Assume that $|1-2x_1x_3|=q^l<q^{2r}$, or $l\leq 2r-1$.
Since $|x_2|\leq q^r$, the condition (\ref{7.16}) cannot be satisfied
for $l=2r-1$. Hence $l\leq 2r-2$. Furthermore, (\ref{7.16}) implies that
$|1-2x_1x_3|q^{-r}<|x_2|<q^r$. Hence, the set (\ref{7.18}) is
\EQ\label{7.18a}
|x_1|\leq q^r,\ |x_3|\leq q^r,\ \frac{|1-2x_1x_3|}{q^{r-1}}\leq |x_2|
\leq q^ {r-1}.
\EN
Set $|x_2|=q^{r_1}$. By Lemma 3.3, we have
$$
\int\limits_{q^{-r}\leq |\al|\leq q^{-r_1}}|\al|^3\psi(\al x_2)\chi(\al)
d^\times\al=\sum_{k=-r}^{-r_1}(q^3X)^k
=\frac{q^3X}{q^3X-1}(q^3X)^{-r_1}-\frac{(q^3X)^{-r}}{q^3X-1}.
$$
The other integral over $\al$ in (\ref{7.17}) is
$$
\int\limits_{|\al|=q|x_2|^{-1}} |\al|^3\psi(\al x_2)\chi(\al)
d^\times\al=-\cir^{-1}q^2Xq^{-3r_1}X^{-r_1}=-\frac{1}{q-1}q^3X(q^3X)^{-r_1}.
$$
Summing up the last two integrals and simplifying the result,
the integration in (\ref{7.17}) over $x_2$ is
$$
\int\limits_{q^{l-r+1}\leq |x_2|\leq q^{r-1}}
\left[\frac{q^4X(1-q^2X)}{(q^3X-1)(q-1)}
(q^3X)^{-r_1}-\frac{(q^3X)^{-r}}{q^3X-1}\right]dx_2.
$$
This integral is equal to
$$
\cir\frac{q^4X(1-q^2X)}{(q^3X-1)(q-1)}
\sum_{r_1=l-r+1}^{r-1}\biggl(\frac{1}{q^2X}\biggr)^{r_1}
-\cir\frac{(q^3X)^{-r}}{q^3X-1}\sum_{r_1=l-r+1}^{r-1}q^{r_1}.
$$
Once simplified, it is equal to
$$
\frac{1}{q}(q^3X+1)(q^2X)^{-r}-\frac{q(q^2X)^{r+1}}{q^3X-1}(q^2X)^{-l}
+\frac{(q^4X)^{-r}}{q^3X-1}q^l.
$$
Hence, the integral (\ref{7.17}) is
\EQ\label{7.19a}
\frac{1}{q}(q^3X+1)(q^2X)^{-r}\sum_{l\leq 2r-2}\int\!\!\int_{|x_1|\leq q^r\!,\,
|x_3|\leq q^r\!,\,|1-2x_1x_3|=q^{l}}dx_1dx_3
\EN
\EQ\label{7.19b}
-\sum_{l\leq 2r-2}\left[\frac{q(q^2X)^{r+1}}{q^3X-1}(q^2X)^{-l}-
\frac{(q^4X)^{-r}}{q^3X-1}q^l\right]\int\!\!\int_{|x_1|\leq q^r\!,\,
|x_3|\leq q^r\!,\,|1-2x_1x_3|=q^l}dx_1dx_3.
\EN
The sums (\ref{7.15b}) and (\ref{7.19b}) cancel each other. Thus,
the sum of the contributions of Case 1b(ii) and Case 2a is obtained on
adding (\ref{7.15a}) and (\ref{7.19a}). It is
$$
\frac{1}{q}(q^3X+1)(q^2X)^{-r}\sum_{l\leq 2r-1}\int\!\!\int_{|x_1|\leq q^r,\,
|x_3|\leq q^r,\,|1-2x_1x_3|=q^{l}}dx_1dx_3.
$$
The sum of the integrals in this expression is the integral
$$
\int\!\!\int_{|x_1|\leq q^r\!,\,|x_3|\leq q^r\!,\,|1-2x_1x_3|
\leq q^{2r-1}}dx_1dx_3
$$
$$
=\int_{|x_1|\leq q^{r-1}}dx_1\int_{|x_3|\leq q^r}dx_3+
\int_{|x_1|=q^r}dx_1\int_{|x_3|\leq q^{r-1}}dx_3=q^{2r-1}(2-1/q).
$$
In conclusion, the contribution from Case 1b(ii) and Case 2a is
\EQ\label{7.20}
\frac{1}{q}(q^3X+1)(q^2X)^{-r}q^{2r-1}\biggl(2-\frac{1}{q}\biggr).
\EN
\indent
{\em Case 2b.} Now $|1-2x_1x_3|=q^{2r}$. Since $|x_i|\leq q^r$, this implies
that $|x_1|=|x_3|=q^r$. Futhermore, to satisfy (\ref{7.16}), we must have
that $|x_2|=q^r$.  Since in this case the inequality
$q^{-r}\leq |\al|\leq |x_2|^{-1}$ is equivalent to  $|\al|=q^{-r}$,
and $q|x_2|^{-1}=q^{1-r}$,
the integral (\ref{7.17}) is equal to
\EQ\label{7.21}
\int\!\!\int\!\!\int\left[\int_{|\al|=q^{-r}}\!\!
+\int_{|\al|=q^{1-r}}\right]|\al|^3\psi(\al x_2)\chi(\al)d^\times\al
dx_1dx_2dx_3.
\EN
where $x_1$, $x_2$ and $x_3$ range over the set
\EQ\label{7.22}
|x_1|=q^r,\ |x_2|=q^r,\ |x_3|=q^r,\
\left|\frac{1-2x_1x_3}{x_2}-x_2\right|\leq q^{r-1}.
\EN
Since $|x_2|=q^r$, using Lemma 3.3, we have
\EQ\label{7.23}
\int_{|\al|=q^{-r}}|\al|^3\psi(\al x_2)\chi(\al)d^\times\al=(q^3X)^{-r},
\EN
and
\EQ\label{7.24}
\int_{|\al|=q^{1-r}}|\al|^3\psi(\al x_2)\chi(\al)d^\times\al
=-\frac{(q^3X)^{1-r}}{q-1}.
\EN
The volume of the set (\ref{7.22}) can be computed as follows.
The inequality (\ref{7.16}) is equivalent to
$$
1-2x_1x_3-x_2^2=\pii^{1-r}tx_2,\ \mbox{  where }\ |t|\leq 1.
$$
Thus, the set (\ref{7.22}) can be described as
$$
|x_1|=q^r,\ |x_2|=q^r,\ x_3=\frac{1-x_2^2}{2x_1}-
\frac{x_2\pii^{1-r}}{2x_1}t,\ |t|\leq 1.
$$
Note that $|x_3|=q^r$ and $dx_3=q^{r-1}dt$. The volume of
(\ref{7.22}) is
$$
q^{r-1}\int_{|x_1|=q^r}dx_1\int_{|x_2|=q^r}dx_2\int_{|t|\leq 1}dt=
\cir^2q^{3r-1}.
$$
Multiplying this by the sum of (\ref{7.23}) and (\ref{7.24}), the
integral (\ref{7.21}) is
$$
\left[(q^3X)^{-r}-\frac{(q^3X)^{1-r}}{q-1}\right]\cir^2q^{3r-1}.
$$
Hence, the contribution to (\ref{7.7}) of Case 2b is
\EQ\label{7.25}
\frac{1}{q}\cir^2X^{-r}-(q-1)X^{1-r}.
\EN
\newline\indent
We have considered all possible cases. Finally, the answer (the integral
(\ref{7.7})) is obtained
on adding (\ref{7.11}), (\ref{7.12}), (\ref{7.20}) and (\ref{7.25}).
\hspace*{\fill}$\Box$

\noindent{\bf Proposition 3.6.} {\em In the non-split case, the Fourier
transform $\widehat{\Psi}_r(\chi)$ of $\Psi(\al,\Phi_r)$ is equal to
$$
X^{-r}-qX^{1-r}
$$
if the character $\chi$ is unramified, and to $0$ if $\chi$ is ramified.}

\noindent{\em Proof.}
Recalling the definitions of $a$, $n$, $\gamma_0$ and $v_0$,
we have
$$
\Psi(\al,\Phi_r)=\int\!\!\!\int\!\!\!\int_{F^3}\Phi_r(2\al\th-x_1+z/\al,x_1/\al,
x_2/\al,1+x_3/\al,1/\al)\psi(x_1+2\th x_3)dx_1dx_2dx_3.
$$
Hence, by definition the integral $\widehat{\Psi}_r(\chi)$ is
$$
\int_{F^\times}\int_{F^3}\Phi_r(2\al\th-x_1+z/\al,x_1/\al,
x_2/\al,1+x_3/\al,1/\al)\psi(x_1+2\th x_3)dx_1dx_2dx_3\chi(\al)d^\times\al.
$$
Recall that $\Phi_r(x_1,x_2,x_3,x_4,x_5)$ is $1$ if $|x_i|\leq q^r$
$(i=1,...,5)$ and is zero otherwise. Thus, the integral above
becomes
$$
\int\!\!\!\int\!\!\!\int\!\!\!\int
\chi(\al)\psi(x_1+2\th x_3)dx_1dx_2dx_3d^\times\al,
$$
over the set defined by
$$
\left|2\al\th-x_1+\frac{z}{\al}\right|\leq q^r,\
\left|\frac{x_1}{\al}\right|\leq q^r\!,\
\left|\frac{x_2}{\al}\right|\leq q^r\!,\
\left|1+\frac{x_3}{\al}\right|\leq q^r\!,\
\left|\frac{1}{\al}\right|\leq q^r\!,
$$
where $z=-x_1x_3-x_2^2/2$.
After the change of variables $x_i\mapsto \al x_i$ followed by $x_3\mapsto 1+x_3$
we arrive at the integral
\EQ\label{B1.1}
\int\!\!\!\int\!\!\!\int\!\!\!\int
|\al|^3\psi(\al (x_1+2\th(x_3-1))\chi(\al)d^\times\al dx_1dx_2dx_3
\EN
over the set given by
\EQ\label{B1.2}
|x_1|\leq q^r,\ |x_2|\leq q^r,\ |x_3|\leq q^r,\
q^{-r}\leq |\al|\leq q^r|2+\th z|^{-1},
\EN
where $z=-x_1x_3-x_2^2/2$. By Lemma 3.3
\EQ\label{B1.3}
\int_{|\al|=q^k}|\al|^3\psi(\al x)\chi(\al)d^\times\al=\left\{
\begin{array}{ll} 0, & \mbox{ if }\ q^k\geq q^2|x|^{-1}, \\
-(q-1)^{-1}(q^3X)^k, & \mbox{ if }\ q^k=q|x|^{-1}, \\
(q^3X)^k, & \mbox{ if }\ q^k\leq |x|^{-1}.
\end{array}\right.
\EN
In order to use this,
we will split (\ref{B1.2}) into two subdomains: according to whether
$q^r/|2+\th z|$ is $\leq$ or $>$ than $1/|x_1+2\th(x_3-1)|$.

{\em Case 1.} We have $q^r/|2+\th z|\leq 1/|x_1+2\th(x_3-1)|$.
Using (\ref{B1.3}), the integration over $\al$ in (\ref{B1.1}) is
performed when $q^{-r}\leq |\al|\leq q^r|2+\th z|^{-1}$. The integral
(\ref{B1.1}) is equal to
$$
\int\!\!\!\int\!\!\!\int\int_{q^{-r}\leq |\al|\leq q^r|2+\th z|^{-1}}
|\al|^3\psi(\al (x_1+2\th(x_3-1))\chi(\al)d^\times\al dx_1dx_2dx_3,
$$
where the triple integral is taken over the set defined by
$$
|x_1|\leq q^r,\ |x_2|\leq q^r,\ |x_3|\leq q^r,\
|x_1+2\th(x_3-1)|\leq |2+\th z|q^{-r}.
$$
Define $l$ by $|2+\th z|=q^l$. The integral above can be writen as
\EQ\label{B1.4}
\sum_{l\leq 2r}q^{l-r}\int\!\!\!\int\!\!\!\int
\int_{q^{-r}\leq |\al|\leq q^{r-l}}
|\al|^3\psi(\al(x_1+2\th(x_3-1))\chi(\al)d^\times\al dtdx_2dx_3,
\EN
where (for each $l$) $x_1$, $x_2$ and $x_3$ range over the set
\EQ\label{B1.5}
|x_1|\leq q^r,\ |x_2|\leq q^r,\ |x_3|\leq q^r,\
|2+\th z|=q^l,\ |x_1+2\th(x_3-1)|\leq q^{l-r}.
\EN
Applying (\ref{B1.3}), we obtain
$$
\int_{q^{-r}\leq |\al|\leq q^{r-l}}
|\al|^3\psi(\al(x_1+2\th(x_3-1)))\chi(\al)d^\times\al=
\sum_{k=-r}^{r-l}(q^3X)^k
$$
\EQ\label{B1.6}
=\frac{(q^3X)^{r-l+1}-(q^3X)^{-r}}{q^3X-1}.
\EN
Making the change of variables in (\ref{B1.5}),
$x_1=2\th(1-x_3)+\pii^{r-l}t$, where $|t|\leq 1$, the sum (\ref{B1.4})
is equal to
\EQ\label{B1.7}
\sum_{l\leq 2r}\mbox{vol}(V_1(l))q^{l-r}\frac{(q^3X)^{r-l+1}-(q^3X)^{-r}}
{q^3X-1},
\EN
where $V_1(l)$ is the set defined by
\EQ\label{B1.8}
|t|\leq 1,\ |x_2|\leq q^r,\ |x_3|\leq q^r,\
|4-\th^2-2\th\pii^{r-l}tx_3+(2\th x_3-\th)^2-\th x_2^2|=q^l.
\EN
Note that since $|x_3|\leq q^r$ and $|t|\leq 1$, we have
$|2\th\pii^{r-l}tx_3|\leq q^l$. We distinguish between the following
subcases.
\np
{\em Case 1a.} Assume that $|2\th\pii^{r-l}tx_3|=q^l$. It follows that
$|t|=1$ and $|x_3|=q^r$. This subset of (\ref{B1.8}) is given by
\EQ\label{B1.9}
|t|=1,\ |x_2|\leq q^r,\ |x_3|=q^r,\
|4-\th^2-2\th\pii^{r-l}tx_3+(2\th x_3-\th)^2-\th x_2^2|=q^l.
\EN
Since $\th$ is a non-square element and $|x_2|\leq q^r$,
$|(2\th x_3-\th)^2-\th x_2^2|=q^{2r}$. Thus the only $l$ when the set
(\ref{B1.9}) is non-empty is when $l=2r$. Once simplified, it is defined by
$$
|t|=1,\ |x_2|\leq q^r,\ |x_3|=q^r,\ |2\pii^{-r}tx_3+x_2^2-\th(2x_3)^2|=q^{2r}.
$$
Alternatively, once $x_2$ and $x_3$ are fixed, $t$ can be any element
with $|t|=1$ which does not belong to $(\th(2x_3)^2-x_2^2)\pii^r/(2x_3)+\pii R$.
Thus, the volume of (\ref{B1.9}) is equal to
\EQ\label{B1.10}
q^{2r}\cir\biggl(1-\frac{2}{q}\biggr).
\EN
\np
{\em Case 1b.} Assume that $|2\th\pii^{r-l}tx_3|<q^l$.
This subset of (\ref{B1.8}) is given by
\EQ\label{B1.11}
|t|\leq 1,\ |x_2|\leq q^r,\ |x_3|=q^r,\ |tx_3|< q^r,\
|4-\th^2+(2\th x_3-\th)^2-\th x_2^2|=q^l.
\EN
Removing the third inequality, we enlarge the set (\ref{B1.11}) by
\EQ\label{B1.12}
|t|\leq 1,\ |x_2|\leq q^r,\ |x_3|=q^r,\ |tx_3|=q^r,\
|4-\th^2+(2\th x_3-\th)^2-\th x_2^2|=q^l.
\EN
Note this subset is non-empty only when $l=2r$, in which case
its volume is
\EQ\label{B1.13}
\int_{|t|=1|}dt\int_{|x_2|\leq q^r}dx_2\int_{|x_3|=q^r}dx_3
=\cir^2q^{2r}.
\EN
Thus, when $l<2r$, the set (\ref{B1.11}) is given by
\EQ\label{B1.14}
|t|\leq 1,\ |x_2|\leq q^r,\ |x_3|=q^r,\ |tx_3|< q^r,\
|4-\th^2+(2\th x_3-\th)^2-\th x_2^2|=q^l,
\EN
and, when $l=2r$ it is the difference of (\ref{B1.13}) and (\ref{B1.12}).

We obtained that when $l<2r$, $\mbox{vol}(V_1(l))=W_l-W_{l-1}$, and
when $l=2r$,
$$
\mbox{vol}(V_1(2r))=W_{2r}-W_{2r-1}+\mbox{(\ref{B1.10})}
-\mbox{(\ref{B1.13})}.
$$
Note that when $l=2r$, (\ref{B1.6}) is equal to $(q^3X)^{-r}$ and
(\ref{B1.10})$-$(\ref{B1.13}) is
$$
q^{2r}\cir\biggl(1-\frac{2}{q}\biggr)-\cir^2q^{2r}=-\cir\frac{1}{q}q^{2r}.
$$
The integral (\ref{B1.1}) over the subset of Case 1 is equal to
\EQ\label{B1.15}
\sum_{l\leq 2r}(W_l-W_{l-1})q^{l-r}\frac{(q^3X)^{r-l+1}-(q^3X)^{-r}}
{q^3X-1}-\cir\frac{1}{q}X^{-r},
\EN
where $W_l$ is the volume of the set defined by
$$
|x_2|\leq q^r,\ |x_3|=q^r,\
|4-\th^2+(2\th x_3-\th)^2-\th x_2^2|\leq q^l.
$$
\np
{\em Case 2.} We have  $q^r/|2+\th z|> 1/|x_1+2\th(x_3-1)|$.
Using (\ref{B1.3}), the integration over $\al$ in (\ref{B1.1}) is
performed when $q^{-r}\leq |\al|\leq q|x_1+2\th(x_3-1)|^{-1}$.
The integral (\ref{B1.1}) is equal to
$$
\int\!\!\!\int\!\!\!\int\int_{q^{-r}\leq |\al|\leq q|x_1+2\th(x_3-1)|^{-1}}
|\al|^3\psi(\al (x_1+2\th(x_3-1))\chi(\al)d^\times\al dx_1dx_2dx_3,
$$
where the triple integral is taken over the set defined by
$$
|x_1|\leq q^r,\ |x_2|\leq q^r,\ |x_3|\leq q^r,\
|x_1+2\th(x_3-1)|> |2+\th z|q^{-r}.
$$
Define $l$ by $|x_1+2\th(x_3-1)|=q^l$. The integral above can be writen as
\EQ\label{B1.16}
\sum_{l\leq r}q^{l-r}\int\!\!\!\int\!\!\!\int
\int_{q^{-r}\leq |\al|\leq q|x_1+2\th(x_3-1)|^{-1}}
|\al|^3\psi(\al(x_1+2\th(x_3-1))\chi(\al)d^\times\al dtdx_2dx_3,
\EN
where (for each $l$) $x_1$, $x_2$ and $x_3$ range over the set
\EQ\label{B1.17}
|x_1|\leq q^r,\ |x_2|\leq q^r,\ |x_3|\leq q^r,\
|x_1+2\th(x_3-1)|= q^,\ |2+\th z|<q^{l+r}.
\EN
Applying (\ref{B1.3}), we split the integral over $\al$ into two integrals
$$
\biggl[\int_{q^{-r}\leq |\al|\leq q^{-l}}+\int_{|\al|=q^{1-l}}\biggr]
|\al|^3\psi(\al(x_1+2\th(x_3-1)))\chi(\al)d^\times\al=
\sum_{k=-r}^{-l}(q^3X)^k-\frac{(q^3X)^{1-l}}{q-1}
$$
\EQ\label{B1.18}
=\frac{(q^3X)^{1-l}-(q^3X)^{-r}}{q^3X-1}-\frac{(q^3X)^{1-l}}{q-1}.
\EN
Making the change of variables in (\ref{B1.17}),
$x_1=2\th(1-x_3)+\pii^{-l}\ep$, where $|\ep|=1$, the sum (\ref{B1.16})
is equal to
\EQ\label{B1.19}
\sum_{l\leq r}\mbox{vol}(V_2(l))q^l\biggl[\frac{(q^3X)^{1-l}-(q^3X)^{-r}}
{q^3X-1}-\frac{(q^3X)^{1-l}}{q-1}\biggr],
\EN
where $V_2(l)$ is the set defined by
\EQ\label{B1.20}
|\ep|=1,\ |x_2|\leq q^r,\ |x_3|\leq q^r,\
|4-\th^2-2\th\pii^{r-l}\ep x_3+(2\th x_3-\th)^2-\th x_2^2|=q^l.
\EN
Note that since $|x_3|\leq q^r$ and $|\ep|=1$, we have
$|2\th\pii^{-l}\ep x_3|\leq q^{r+l}$. We distinguish between the following
subcases.
\np
{\em Case 2a.} Assume that $|2\th\pii^{-l}\ep x_3|=q^{l+r}$. It implies that
$|x_3|=q^r$. Following the same argument as in Case 1a, we conclude that
with this assumption, the only non-empty subset of (\ref{B1.20}) is when
$l=r$. It is defined by
\EQ\label{B1.21}
|\ep|=1,\ |x_2|\leq q^r,\ |x_3|=q^r,\
|2\pii^{-r}\ep x_3+x_2^2-\th(2x_3)^2|\leq q^{2r-1}.
\EN
Alternatively, once $x_2$ and $x_3$ are fixed, $\ep$ should be in
$-(x_2^2-\th(2x_3)^2)\pii^r/(2x_3)+\pii R$. Note that $|\ep|=1$.
The volume of (\ref{B1.21}) is
$$
q^{2r}\cir\frac{2}{q}.
$$
When $l=r$ (\ref{B1.18}) is equal to $(q^3X)^{-r}-(q^3X)^{1-r}/(q-1)$.
The contribution of this case to (\ref{B1.19}) is
\EQ\label{B1.22}
q^{2r}\cir\frac{2}{q}q^r\biggl[(q^3X)^{-r}-\frac{(q^3X)^{1-r}}{q-1}\biggr]
\EN
\indent
{\em Case 2b.} Assume that $|2\th\pii^{-l}\ep x_3|\leq q^{l+r-1}$.
Thus this subset of (\ref{B1.20}) is given by
\EQ\label{B1.23}
|\ep|= 1,\ |x_2|\leq q^r,\ |x_3|=q^r,\
|4-\th^2+(2\th x_3-\th)^2-\th x_2^2|\leq q^{l+r-1}.
\EN
Note that the volume of this set is $\cir W_{l+r-1}$.

Combining these two subcases,
the integral (\ref{B1.1}) over the subset of Case 2 is equal to
\EQ\label{B1.24}
\sum_{l\leq r}W_{l+r-1}q^l\cir\biggl[\frac{(q^3X)^{1-l}-(q^3X)^{-r}}
{q^3X-1}-\frac{(q^3X)^{1-l}}{q-1}\biggr]+\frac{1}{q}\cir X^{-r}-qX^{1-r}.
\EN
The answer is obtained on adding (\ref{B1.15}) to (\ref{B1.24}). Fix any
$k<2r$. To find the coefficient of $W_k$ in (\ref{B1.15}), we consider
the terms when $l=k$ and $l=k+1$. This coefficient is equal to
\EQ\label{B1.25}
q^{k-r}\frac{(q^3X)^{r-k+1}-(q^3X)^{-r}}{q^3X-1}-
q^{k+1-r}\frac{(q^3X)^{r-k}-(q^3X)^{-r}}{q^3X-1}.
\EN
Similarly, to find the coefficient of $W_k$ in (\ref{B1.24}), we
consider the term with $l=1+k-r$. The coefficient is equal to
$$
\cir q^{1+k-r}\left[\frac{(q^3X)^{r-k}-(q^3X)^{-r}}{q^3X-1}-
\frac{(q^3X)^{r-k}}{q-1}\right]
$$
\EQ\label{B1.26}
=q^{k+1-r}\frac{(q^3X)^{r-k}-(q^3X)^{-r}}{q^3X-1}-
q^{k-r}\frac{(q^3X)^{r-k}-(q^3X)^{-r}}{q^3X-1}-q^{k-r}(q^3X)^{r-k}.
\EN
The second term of (\ref{B1.25}) cancels the first one of (\ref{B1.26}).
Thus, their sum is zero. Futher, note that
$$
W_{2r}=\int_{|x_2|\leq q^r}dx_2\int_{|x_3|=q^r}dx_3=\cir q^{2r}.
$$
Thus the sum of (\ref{B1.15}) and (\ref{B1.24}) is equal to
$$
W_{2r}q^r\frac{(q^3X)^{1-r}-(q^3X)^{-r}}{q^3X-1}-\cir\frac{1}{q}X^{-r}
+\cir\frac{1}{q}X^{-r}-qX^{1-r}=X^{-r}-qX^{1-r}.
$$
The Proposition is proved.
\hspace*{\fill}$\Box$

\noindent
{\bf Theorem.} {\em Corresponding $f$ and $f'$ are matching. }

\noindent
{\em Proof.} Indeed as we have seen in Section I.0,
to prove that corresponding functions are
matching (i.e. $\Psi(\al,\phi_f)=\psi(\al)|\al|\Psi'(\al^{-1},\phi'_{f'})$)
it is enough to show that (for $r\geq 0$)
\EQ\label{7:26}
\Psi(\al,\Phi_r)=\psi(\al)|\al|
\Psi'(\al^{-1},(-1)^rq^r(\phi'_r\pm \phi'_{r-1})).
\EN
\np
Comparing Proposition 3.5 in the split case, and Proposition 3.6 in the
non-split case, with Proposition 3.1, we have (for $r\geq 1$)
$$
\int_{F^\times}\Psi(\al,\Phi_r)\chi(\al)d^\times\al=
\int_{F^\times}\psi(\al)|\al|\Psi'(\al^{-1},(-1)^rq^r
(\phi'_r\pm \phi'_{r-1}))
\chi(\al)d^\times\al,
$$
where $\chi$ is any complex valued character of $F^\times$.
If $\chi$ is ramified both integrals are equal to $0$.
Fourier inversion formula now implies (\ref{7:26}) when $r\geq 1$.
When $r=0$, the formula (\ref{7:26}) follows from the
unit element case, treated in [FM].
\hspace*{\fill}$\Box$

\newpage
\begin{center}{\bf References}\end{center}

[BZ] I. Bernstein, A. Zelevinskii, Representations of the group $GL(n,F)$
where $F$ is a non-archimedean local field, {\em Uspekhi Mat. Nauk},
31 (1976), 5-70.

[F] Y. Flicker, Twisted Tensors and Euler Products,
{\em Bull. Soc. math. France,} 116 (1988), 295-313.

[FM] Y. Flicker, J. G. M. Mars, Cusp forms on $GSp(4)$ with
$SO(4)$-periods, {\em Preprint.}

[J] H. Jacquet, Relative Kloosterman Integrals for GL(3): II,
{\em Canad. J. Math.} 44 (1992), 1220-1240.

[L] R.P. Langlands, Automorphic representations, Shimura varieties,
and motives, {\em Proc. Sympos. Pure Math.} 33 II (1979), 205-246.

[M] Z. Mao, Relative Kloosterman Integrals for GL(3): III,
{\em Canad. J. Math.} 45 (1993), 1211-1230.

[MS] A. Murase, T. Sugano, Shintani function and its application to automorphic
$L$-functions for classical groups, {\em Math. Ann.} 299 (1994), 17-56.

[O] T. Oda, On modular forms associated with indefinite quadratic forms of
signature $(2,n-2)$, {\em Math. Ann.} 231 (1977), 97-144.

[PS] I. Piatetski-Shapiro, On the Saito-Kurokawa lifting, {\em Invent. Math.}
71 (1983), 309-338; 76 (1984), 75-76.

[R] S. Rallis, On a relation between SL(2) cusp forms and automorphic
forms on orthogonal groups, {\em Proc. Sympos. Pure  Math.} 33 I (1979), 297-314.

[RS] S. Rallis, G. Schiffmann, On a relation between SL(2) cusp forms
and cusp forms on tube domains associated to orthogonal groups,
{\em Trans. AMS} 263 (1981), 1-58.

[Sh] T. Shintani, On an explicit formula for class 1 ``Whittaker functions''
on $GL_n$ over $p$-adic fields, {\em Proc. Japan Acad.} 52 (1976), 180-182.

[T]  J. Tits, Reductive groups over local fields,
{\em Proc. of Sympos. Pure Math.} 33 I (1979), 29-69.

\end{document}